\documentclass[11pt]{article}
\usepackage[a4paper,margin=25mm]{geometry}
\usepackage{amsmath,amssymb,bm}
\usepackage{booktabs}
\usepackage{tabularx}
\usepackage{longtable}
\usepackage{graphicx}
\usepackage{adjustbox}
\usepackage{subcaption}
\usepackage{placeins}
\usepackage{tikz}
\usetikzlibrary{arrows.meta,backgrounds,calc,fit,positioning}
\usepackage{microtype}
\usepackage{siunitx}
\usepackage{xcolor}
\usepackage{hyperref}
\hypersetup{colorlinks=true,linkcolor=blue,citecolor=blue,urlcolor=blue}
\DeclareSIUnit{\mebibyte}{MiB}
\newcommand{\dd}{\mathop{}\!\mathrm{d}}
\newcommand{\R}{\mathbb{R}}
\newcommand{\norm}[1]{\left\lVert#1\right\rVert}
\newcommand{\ResultTableFont}{\small}
\newcommand{\VerificationSchematicWidth}{0.98\linewidth}
\newcommand{\FullFigureWidth}{0.92\linewidth}

\newcommand{\WedgeNormalizedInput}{\ensuremath{x=(\alpha/\pi-0.55)/0.45}}

\newcommand{\WedgeParameterBounds}{\ensuremath{[-8,8]}}
\newcommand{\WedgeRegularization}{\ensuremath{10^{-7}}}

\newcommand{\WedgeTrainAngles}{\ensuremath{(0.58,0.70,0.82,0.94,1.00)}}
\newcommand{\WedgeTestAngles}{\ensuremath{(0.64,0.76,0.88,0.97)}}
\newcommand{\WedgeTrainableCount}{\ensuremath{4}}

\newcommand{\NonlinearTrainingControls}{\ensuremath{(0.35,0.55,0.80,1.05,1.25)}}
\newcommand{\NonlinearTestControls}{\ensuremath{(0.30,1.40)}}

\newcommand{\LearnedQ}{\ensuremath{2.000007}}
\newcommand{\LearnedWeightOne}{\ensuremath{0}}
\newcommand{\LearnedWeightTwo}{\ensuremath{0}}

\newcommand{\VectorDOF}{\ensuremath{156}}

\newcommand{\VectorYoungModulus}{\ensuremath{1}}
\newcommand{\VectorPoissonRatio}{\ensuremath{0.3}}
\newcommand{\VectorBenchmarkRadius}{\ensuremath{1}}
\newcommand{\VectorTargetKI}{\ensuremath{1.25}}
\newcommand{\VectorTargetKII}{\ensuremath{-0.45}}
\newcommand{\VectorTargetTStress}{\ensuremath{0.08}}

\newcommand{\VectorLearnedSIFCorrelationError}{\ensuremath{9.931\times 10^{-5}}}
\newcommand{\VectorEnergyGain}{\ensuremath{33.42}}
\newcommand{\VectorSIFCorrelationGain}{\ensuremath{8648}}
\newcommand{\FixedLearnedEnergyRatio}{\ensuremath{1}}
\newcommand{\FixedLearnedSIFRatio}{\ensuremath{0.9787}}

\newcommand{\LearnedMinimumEigenvalue}{\ensuremath{0.03931}}
\newcommand{\LearnedMinimumSampledJacobian}{\ensuremath{1.501\times 10^{-7}}}

\newcommand{\WedgeANNWorstQError}{\ensuremath{0.01082}}

\newcommand{\WedgeAffineWorstQError}{\ensuremath{1.667\times 10^{-4}}}

\newcommand{\InteractionTargetKI}{\ensuremath{1.25}}
\newcommand{\InteractionTargetKII}{\ensuremath{-0.45}}

\newcommand{\InteractionContourCount}{\ensuremath{3}}

\newcommand{\InteractionLearnedWorstError}{\ensuremath{1.923\times 10^{-5}}}

\newcommand{\InteractionLearnedKISpread}{\ensuremath{1.48\times 10^{-5}}}
\newcommand{\InteractionLearnedKIISpread}{\ensuremath{2.988\times 10^{-6}}}

\newcommand{\InteractionExactMaxError}{\ensuremath{1.091\times 10^{-14}}}

\newcommand{\NonlinearHalfAngleOverPi}{\ensuremath{0.75}}
\newcommand{\NonlinearBeta}{\ensuremath{3}}
\newcommand{\NonlinearModulusMultiplier}{\ensuremath{4}}
\newcommand{\NonlinearModulusCoefficient}{\ensuremath{12}}
\newcommand{\NonlinearTrainingMinimum}{\ensuremath{0.35}}
\newcommand{\NonlinearTrainingMaximum}{\ensuremath{1.25}}
\newcommand{\NonlinearTestMinimum}{\ensuremath{0.3}}
\newcommand{\NonlinearTestMaximum}{\ensuremath{1.4}}

\newcommand{\NonlinearRepresentativeControl}{\ensuremath{1.4}}
\newcommand{\NonlinearAffineResidual}{\ensuremath{0.1953}}

\newcommand{\NonlinearEnergyGainMinimum}{\ensuremath{4.012}}
\newcommand{\NonlinearFieldGainMinimum}{\ensuremath{13.34}}
\newcommand{\NonlinearAmplitudeGainMinimum}{\ensuremath{1.123}}

\newcommand{\NonlinearDensityIncrementalEnergyGainMinimum}{\ensuremath{4.839}}
\newcommand{\NonlinearDensityOnPowerQEnergyGainMinimum}{\ensuremath{4.732}}

\newcommand{\MapQuadratureEnergyDelta}{\ensuremath{3.681\times 10^{-8}}}
\newcommand{\AssemblyQuadratureEnergyDelta}{\ensuremath{5.738\times 10^{-11}}}

\newcommand{\PrimaryVectorEnergyMean}{\ensuremath{0.005057}}
\newcommand{\PrimaryVectorEnergyStd}{\ensuremath{1.484\times 10^{-8}}}

\newcommand{\JacobianCaseCount}{\ensuremath{50}}

\newcommand{\JacobianCutoff}{\ensuremath{0.001}}
\newcommand{\JacobianMinimumNormalized}{\ensuremath{2.851\times 10^{-13}}}
\newcommand{\JacobianMaximumDerivativeDisagreement}{\ensuremath{0.004221}}
\newcommand{\JacobianDerivativeTolerance}{\ensuremath{0.005}}

\newcommand{\CostTargetEnergy}{\ensuremath{0.1}}

\newcommand{\AdaptiveRawCondition}{\ensuremath{4.302\times 10^{4}}}
\newcommand{\AdaptiveScaledCondition}{\ensuremath{446.7}}

\newcommand{\AdaptiveFirstSolveSeconds}{\ensuremath{1.914}}

\newcommand{\AdaptiveBreakEven}{no finite break-even}

\title{Mechanics-trained neural coordinate mapping for B-spline analysis of crack-tip and corner singularities}
\author{Hyunju Kim\thanks{Corresponding authors:
hjkim@kentech.ac.kr}\\
Department of Energy Engineering\\
Korea Institute of Energy Technology (KENTECH)\\
Naju 58217, Republic of Korea}
\date{}

\begin{document}
\maketitle

\begin{abstract}
Near a crack tip or re-entrant corner, fractional radial powers can have
unbounded derivatives and slow the convergence of high-order splines.  A
singular mapping grades the computational radius so that the pulled-back field
is smoother without changing the physical domain.  Classical maps require the
singular exponent and grading in advance.  Here the radial coordinate is trained from the mechanics problem.  It is the
normalized integral of a positive neural density, which fixes both radial
boundaries and ensures $r'(s)>0$ away from the collapsed tip.  The radial grading exponent and density-correction weights are inferred from
Galerkin energies evaluated at discrete equilibrium, without exponent labels or
exact interior fields.  We test
the mapping in scalar and plane-strain B-spline formulations and compare it
with the identity map, radially graded knot vectors, prescribed power maps,
and an adaptive
enriched B-spline method.  At \VectorDOF{} vector degrees of freedom, the
mechanics-trained map reduces the relative energy-norm error by a factor of
\VectorEnergyGain{} compared with the identity map.  The maximum error in
the mixed-mode stress intensity factors recovered on
\InteractionContourCount{} contours is \InteractionLearnedWorstError{}.  For
the straight crack, the learned density correction vanishes and the prescribed
$r=s^2$ map gives the same improvement.  For a nonlinear Robin family with test parameters outside the
training interval, the density correction remains nonzero and gives
a fixed-$q$ incremental energy gain of at least
\NonlinearDensityIncrementalEnergyGainMinimum{}.  Thus, for the problems
considered here, the neural correction is useful when the required coordinate
is not represented by a prescribed power map.
\end{abstract}


\section*{Nomenclature}


\begingroup
\footnotesize
\renewcommand{\arraystretch}{0.94}
\setlength{\LTpre}{3pt}
\setlength{\LTpost}{3pt}
\begin{longtable}{@{}p{0.18\linewidth}p{0.76\linewidth}@{}}
\toprule
Symbol & Meaning \\
\midrule
\endfirsthead
\toprule
Symbol & Meaning \\
\midrule
\endhead
\bottomrule
\endfoot
$s,r,q,\lambda$ & computational and physical radii, radial grading exponent, and physical singular exponent \\
$\phi=(q,\bm w),\rho_\phi,r_\phi$ & trainable map parameters, positive radial density, and boundary-normalized radial map \\
$\zeta,\vartheta,q_\vartheta$ & problem parameter, conditional-network parameters, and predicted radial grading exponent \\
$\mathcal J_h$ & discrete strain-energy objective evaluated at Galerkin equilibrium \\
$\bm\chi_\phi,J_\phi$ & implemented physical chart and its integration Jacobian \\
$K_I,K_{II}$ & mode-I and mode-II stress intensity factors \\
\end{longtable}
\endgroup

\section{Introduction}
\label{sec:introduction}

Solutions of elliptic and elasticity problems may contain fractional radial
powers near cracks and re-entrant corners.  If
$u(r,\theta)\sim r^\lambda\Psi(\theta)$ with $0<\lambda<1$, where
$\lambda$ is the leading singular exponent and $\Psi$ its angular eigenfunction,
then the radial derivative is unbounded at $r=0$, and the convergence of a standard high-order
approximation is reduced.  The crack-face discontinuity and the radial
singularity require different treatments.  We retain the two faces as distinct
parametric boundaries and use the coordinate transformation solely to
regularize the factor $r^\lambda$.

Williams' expansion identifies the fractional power responsible for the loss
of regularity at cracks and re-entrant corners \cite{Williams1952Corners}.
This observation led to mappings that transfer the singular behavior to a
smoother computational coordinate \cite{OhBabuska1995}.  The spline geometry
and trial spaces of isogeometric analysis \cite{Hughes2005IGA} make this
transfer especially direct: Jeong et al.\ realized it by collapsing a NURBS
edge onto the singular point \cite{Jeong2013Mapping}.
Repeating the control points on one parametric edge collapses that edge to the
singular point, while the remaining control net retains the physical boundary.
More precisely, let $\widehat\Omega=(0,1)_\xi\times(0,1)_\eta$, where $\xi$
is tangential to the collapsed edge and $\eta$ measures distance from it.  For
control points $\bm P_{ij}$, positive weights $w_{ij}$, and univariate B-spline
bases $N_i$ and $M_j$, the geometry map is
\[
 \bm F(\xi,\eta)=
 \frac{\sum_{i,j}N_i(\xi)M_j(\eta)w_{ij}\bm P_{ij}}
      {\sum_{i,j}N_i(\xi)M_j(\eta)w_{ij}}.
\]
The repeated control layers impose $\bm F(\xi,0)=\bm x_0$ and, for a
$q$th-order collapse, make the first $q-1$ radial derivatives vanish.  Hence,
after translating the singular point $\bm x_0$ to the origin,
\[
 \bm F(\xi,\eta)-\bm x_0
 =\eta^q\widetilde{\bm F}(\xi)+o(\eta^q),
 \qquad
 \widetilde{\bm F}(\xi)
 =\lim_{\eta\downarrow0}
 \frac{\bm F(\xi,\eta)-\bm x_0}{\eta^q}.
\]
Here $\bm F:\widehat\Omega\to\Omega$ is the NURBS geometry map,
$\widetilde{\bm F}$ is its first nonzero radial coefficient, and
$o(\eta^q)$ is a remainder $\bm R_F$ satisfying
$\sup_\xi\|\bm R_F(\xi,\eta)\|/\eta^q\to0$ as $\eta\downarrow0$.
Consequently,
\[
 r\bigl(\bm F(\xi,\eta)\bigr)
 =\eta^q\|\widetilde{\bm F}(\xi)\|[1+o(1)],
 \qquad
 \theta\bigl(\bm F(\xi,\eta)\bigr)
 \longrightarrow \arg\widetilde{\bm F}(\xi).
\]
If the leading physical field is
$u(\bm x)=r^\lambda\Psi(\theta)+o(r^\lambda)$, where $\lambda>0$ is its
singular exponent and $\Psi$ its angular eigenfunction, then
\[
 u\circ\bm F
 =\eta^{q\lambda}\widetilde\Psi(\xi)+o(\eta^{q\lambda}),
 \qquad
 \widetilde\Psi(\xi)=
 \|\widetilde{\bm F}(\xi)\|^\lambda
 \Psi\!\left(\arg\widetilde{\bm F}(\xi)\right).
\]
Thus, if $q\lambda=m\in\mathbb N$, the fractional radial factor is transformed
into the polynomial factor $\eta^m$ in the parameter space.  By the chain rule,
$\partial_r u=(\partial_\eta(u\circ\bm F))/(\partial_\eta r)
\sim\lambda r^{\lambda-1}\Psi$, so the physical singularity is retained while
its pullback becomes regular.  For example, the square-root crack-tip term has
$\lambda=1/2$, and $q=2$ gives a leading term that is linear in $\eta$.
The angular control net and the repeated radial control points can be adjusted
independently, so the physical boundary is retained although the Jacobian
vanishes on the collapsed edge.  Matching $q$ to the prescribed exponent then
recovers the expected approximation order \cite{Jeong2013Mapping}.  Enriched
and locally refined formulations pursue the same regularity gain by modifying
the trial space instead of the coordinate
\cite{OhKimJeong2014,Fries2010XFEM,Gu2019XIGA}; robust IGA formulations likewise
treat the singular map as a prescribed part of the discretization
\cite{Jonsson2025SingularMaps}.  In all of these approaches, the relevant
singular structure is supplied before the discrete problem is solved.

Prescribing $q$ is sufficient when the leading exponent and a suitable power
map, such as $r=s^2$ for a square-root field, are known.  It becomes restrictive
for a parameterized problem in which the assembled mechanics must reveal both
the exponent and a non-power-law radial correction.  We therefore determine these
quantities from Galerkin energies evaluated at discrete equilibrium, without field or exponent
labels, while enforcing the radial anchors and monotonicity analytically.

Most neural fracture formulations place the network in the approximation of
the physical field.  Enriched PINNs, for example, build crack-tip functions
into a neural trial space \cite{Gu2023EnrichedPINN,Zhu2025ExtremeFracture},
whereas peridynamic, phase-field, and finite-element neural models use it to
represent fracture evolution or material response
\cite{Ning2023PeridynamicPINN,Dammass2025Hybrid,Pantidis2026IFENN}.  Networks
have also been used to transform domains for geometric normalization and
differentiable PDE solution
\cite{Zhan2025Parameterization,Chen2026JacobiNet,Gasick2023IGINN}.  These two
uses meet at a different question: can the network improve a Galerkin space by
changing its coordinate without replacing its unknown field?  We answer this
question by assigning the network only to the derivative of the singular
coordinate.  The displacement coefficients remain the equilibrium unknowns,
so the assembled stiffness operator and the standard SIF-extraction procedures are
retained.  The coordinate can then be compared with prescribed-power, radially graded-knot,
and enriched-spline discretizations on the same approximation space.

Restricting the trainable component to a one-dimensional coordinate derivative
makes positivity and boundary normalization exact.  It also isolates the
coordinate effect: equilibrium and SIF-extraction functionals remain
unchanged, while admissibility reduces to an analytic radial condition and
Jacobian tests of the resulting chart.

This formulation separates two effects that are easily confused.  Learning
$q$ can recover a classical power map when the leading singularity is already
known to be of power type, whereas learning the positive density can modify the
coordinate when a power map is insufficient.  The straight crack and the
affine wedge isolate the first effect; the nonlinear Robin family, evaluated
outside its training range, isolates the second.  Because every map is trained
through an assembled weak form, these comparisons require neither exact
interior fields nor exponent labels.  We further determine whether the same
coordinate improvement survives independent SIF evaluation,
refinement, quadrature changes, Jacobian tests of the implemented maps, and direct assembly on the collapsed-edge trial space.  The resulting accuracy is finally interpreted together with the
one-time construction cost and the repeated Galerkin-solve cost, using prescribed-power,
radially graded-knot, and adaptively enriched spline discretizations as references.

Section~\ref{sec:method} derives the coordinate and its mechanics-based
training functional; Section~\ref{sec:verification} formulates the verification
problems and output measures; and Section~\ref{sec:results} examines the
resulting accuracy, admissibility, and computational cost.  The conclusions
are summarized in Section~\ref{sec:conclusions}.

\section{Construction of the mechanics-trained singular mapping}
\label{sec:method}

The construction follows from the observation that admissibility is governed
more naturally by the derivative of the radial coordinate than by the
coordinate itself.  A positive derivative, once normalized and integrated,
fixes both radial boundaries and prevents folding.  Its pullback changes the
regularity seen by the B-spline basis, and the Galerkin energy evaluated at the discrete equilibrium solution then
provides a mechanics-based signal for selecting the map.  The same argument
leads to a direct parameterization for one problem and a conditional rule for
a family of problems.  The nonlinear Robin family is introduced at the end of
the section because it distinguishes a genuine density correction from the
selection of a radial grading exponent alone.

\subsection{Construction of a positive radial coordinate}

Let $s\in[0,1]$ denote the computational radius.  We seek
$r_\phi:[0,1]\rightarrow[0,1]$ such that $r_\phi(0)=0$, $r_\phi(1)=1$, and
$r_\phi'(s)>0$ for $s>0$.  Because penalty terms cannot guarantee these
conditions during optimization, we parameterize $r_\phi'$ and recover the
coordinate by normalized integration:
\begin{align}
 \rho_\phi(s)
 &=s^{q-1}\exp\!\left[
   \sum_j w_j\tanh\!\bigl(\beta_j(s-c_j)\bigr)\right],
 \label{eq:density}\\
 r_\phi(s)
 &=\frac{\displaystyle\int_0^s\rho_\phi(\tau)\dd\tau}
         {\displaystyle\int_0^1\rho_\phi(\tau)\dd\tau}.
 \label{eq:radialmap}
\end{align}
In Eqs.~\eqref{eq:density} and \eqref{eq:radialmap}, the exponent $q$ controls
the leading radial grading, $w_j$ are the trainable density-correction weights,
and the fixed $c_j$ and $\beta_j$ are the centers and slopes of the smooth
features.  Since the
exponential factor is positive and $q>0$, the density is positive away from the
prescribed collapsed edge.  The fundamental theorem of calculus then gives the
boundary anchors $r_\phi(0)=0$, $r_\phi(1)=1$, and the derivative
\begin{equation}
 r_\phi'(s)=\frac{\rho_\phi(s)}
 {\displaystyle\int_0^1\rho_\phi(\tau)\dd\tau}>0
 \qquad\text{for }s>0.
 \label{eq:certificate}
\end{equation}
Equation~\eqref{eq:certificate} proves injectivity of the one-dimensional
radial coordinate.  The corresponding multidimensional chart is examined
separately through its Jacobian.

If a leading field behaves as $u\sim r^\lambda\Psi(\theta)$, the pure power
specialization gives $u\circ r_\phi\sim s^{q\lambda}\Psi(\theta)$.  The smallest integer grading exponent that maps the homogeneous square-root
displacement to a linear dependence on $s$ is $q=2$.  This known value is
retained as a prescribed-power baseline; it is never
used as a training label.

\subsection{Galerkin formulation and training functional}

Once $r_\phi$ is fixed, the approximation space is the pullback of tensor
B-splines through the physical chart.  With $a$ denoting the angular coordinate,
we use
\begin{equation}
 V_h(r_\phi)=
 \left\{
 v_h:\ v_h(r_\phi(s),a)=\sum_A N_A(s,a)c_A,\quad c_A\in\R
 \right\}.
 \label{eq:trialspace}
\end{equation}
The space in Eq.~\eqref{eq:trialspace} consists of physical functions whose
pullbacks are tensor-product splines.  Hence the coordinate changes the radial
regularity resolved by the spline basis, whereas $c_A$ remain the usual
Galerkin coefficients.  For a vector field, the corresponding discrete system
is
\begin{equation}
 \bm u_h(s,a)=\sum_{A}N_A(s,a)\bm d_A,
 \qquad
 K_{ff}(\phi)\bm d_f=\bm f_f-K_{fc}(\phi)\bm d_c,
 \label{eq:discrete}
\end{equation}
Here $\bm d_A$ are spline coefficients; subscripts $f$ and $c$ identify free
and prescribed blocks, respectively, so $K_{ff}$ is the free--free stiffness
block, $K_{fc}$ is the coupling block, and $\bm f_f$ is the free load vector.
The constitutive law is homogeneous isotropic plane strain. The crack faces
occupy distinct angular limits even though their physical coordinates meet at the tip.  Their
separation is therefore supplied by the parameter-domain topology, while
$r_\phi$ regularizes the radial behavior.

The training functional is the discrete Dirichlet energy evaluated at the Galerkin solution,
\begin{equation}
 \mathcal J_h(\phi)=\frac{1}{2}\bm d(\phi)^T K(\phi)\bm d(\phi).
 \label{eq:objective}
\end{equation}
The primary scalar training problem is the unit slit disk
\begin{equation}
 -\Delta u=0,\qquad
 u(1,\theta)=\sin(\theta/2),\qquad
 \partial_\theta u(r,\pm\pi)=0,\qquad u(0)=0.
 \label{eq:scalartraining}
\end{equation}
For each $\phi$, Eq.~\eqref{eq:discrete} eliminates the field coefficients and
reduces learning to the scalar functional $\mathcal J_h(\phi)$.  We identify a
direct map from the slit-disk problem and transfer it unchanged to elasticity.
The exact field $u=\sqrt r\sin(\theta/2)$ is used only to evaluate that transfer.
For the wedge family, both the neural and affine conditional rules minimize the
mean of the same reduced energy; their comparison therefore changes the
parameter rule, not the energy objective.

\subsection{Neural parameterization and training procedure}

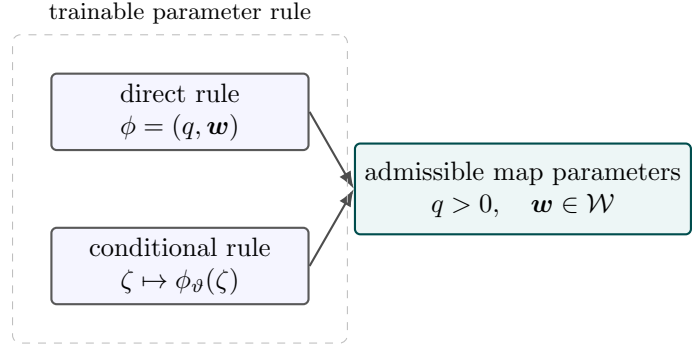
\begin{figure}[t]
\centering
\begin{tikzpicture}[
  font=\small, node distance=10mm and 15mm, >=Latex,
  block/.style={draw=black!65, rounded corners=2pt, thick, align=center,
    minimum height=10mm, minimum width=34mm, fill=blue!4},
  output/.style={draw=teal!60!black, rounded corners=2pt, thick, align=center,
    minimum height=12mm, minimum width=42mm, fill=teal!7},
  flow/.style={->,thick,draw=black!70}]
  \node[block] (direct) {direct rule\\$\phi=(q,\bm w)$};
  \node[block, below=of direct] (conditional)
    {conditional rule\\$\zeta\mapsto\phi_\vartheta(\zeta)$};
  \node[output, right=23mm of $(direct)!0.5!(conditional)$] (parameters)
    {admissible map parameters\\$q>0,\quad \bm w\in\mathcal W$};
  \draw[flow] (direct.east) -- (parameters.west);
  \draw[flow] (conditional.east) -- (parameters.west);
  \node[draw=black!25, dashed, rounded corners=3pt, fit=(direct)(conditional),
    inner sep=5mm, label={[font=\footnotesize]above:trainable parameter rule}] {};
\end{tikzpicture}%
\caption{Direct and conditional parameter rules for map construction.}
\label{fig:parameterrule}
\end{figure}

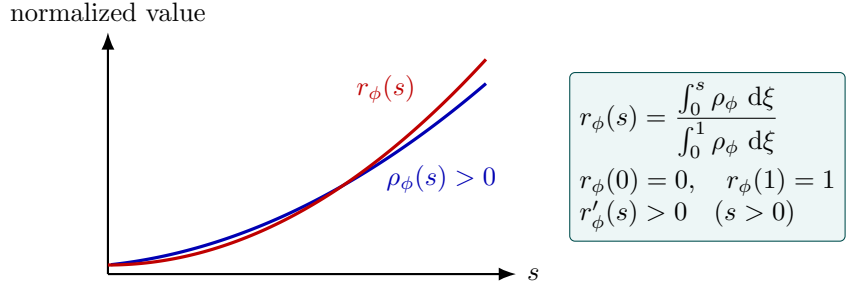
\begin{figure}[t]
\centering
\begin{tikzpicture}[x=1cm,y=1cm,>=Latex,font=\small]
  \draw[->,thick] (0,0) -- (5.4,0) node[right] {$s$};
  \draw[->,thick] (0,0) -- (0,3.2) node[above] {normalized value};
  \draw[blue!70!black,very thick,domain=0.04:5,samples=80]
    plot (\x,{0.12+0.105*\x+0.075*\x*\x});
  \draw[red!75!black,very thick,domain=0:5,samples=80]
    plot (\x,{0.12+2.72*(\x/5)*(\x/5)});
  \node[blue!70!black,anchor=west] at (3.55,1.26) {$\rho_\phi(s)>0$};
  \node[red!75!black,anchor=west] at (3.15,2.45) {$r_\phi(s)$};
  \node[draw=teal!60!black,rounded corners=2pt,fill=teal!7,
    align=left,anchor=west,minimum height=21mm] (normal) at (6.1,1.55)
    {$r_\phi(s)=\dfrac{\int_0^s\rho_\phi\,\dd\xi}
      {\int_0^1\rho_\phi\,\dd\xi}$\\[1mm]
     $r_\phi(0)=0,\quad r_\phi(1)=1$\\
     $r_\phi'(s)>0\quad(s>0)$};
\end{tikzpicture}%
\caption{Positive-density construction of the radial coordinate.}
\label{fig:positivedensity}
\end{figure}

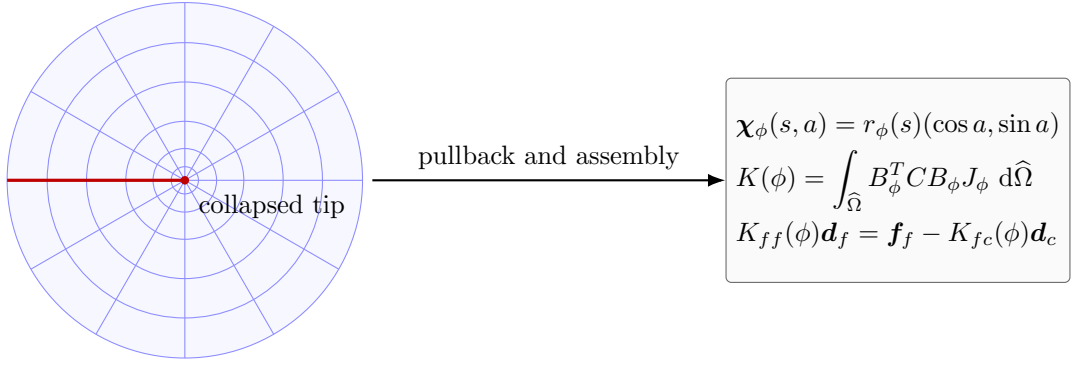
\begin{figure}[t]
\centering
\begin{tikzpicture}[x=1cm,y=1cm,>=Latex,font=\small]
  \coordinate (C) at (0,0);
  \fill[blue!3] (C) circle (2.35);
  \foreach \r in {0.18,0.42,0.78,1.30,1.82,2.35}
    \draw[blue!45] (C) circle (\r);
  \foreach \a in {-150,-120,...,150}
    \draw[blue!45] (C) -- (\a:2.35);
  \draw[very thick,red!75!black] (C) -- (-2.35,0);
  \fill[red!75!black] (C) circle (1.6pt);
  \node[below right=2pt] at (C) {collapsed tip};
  \node[draw=black!60,rounded corners=2pt,fill=gray!4,align=left,
    minimum height=27mm,anchor=west] (assembly) at (7.15,0)
    {$\bm\chi_\phi(s,a)=r_\phi(s)(\cos a,\sin a)$\\[1mm]
     $K(\phi)=\displaystyle\int_{\widehat\Omega}B_\phi^TCB_\phi J_\phi\,\dd\widehat\Omega$\\[1mm]
     $K_{ff}(\phi)\bm d_f=\bm f_f-K_{fc}(\phi)\bm d_c$};
  \draw[->,thick] (2.48,0) -- (assembly.west)
    node[midway,above]{pullback and assembly};
\end{tikzpicture}%
\caption{Mapped Galerkin discretization of the slit disk.}
\label{fig:mappedgalerkin}
\end{figure}

\begin{figure}[t]
\centering
\resizebox{\FullFigureWidth}{!}{%
\begin{tikzpicture}[
  font=\scriptsize, >=Latex,
  train/.style={draw=blue!60!black,rounded corners=2pt,thick,fill=blue!5,
    align=center,text width=22mm,minimum height=18mm,inner sep=2.2mm},
  eval/.style={draw=teal!60!black,rounded corners=2pt,thick,fill=teal!6,
    align=center,text width=18mm,minimum height=18mm,inner sep=2.2mm},
  ref/.style={draw=orange!75!black,rounded corners=2pt,thick,fill=orange!8,
    align=center,text width=20mm,minimum height=18mm,inner sep=2.2mm},
  flow/.style={->,thick,draw=black!72},
  note/.style={font=\bfseries\footnotesize,anchor=south west}]
  \node[train] (parameters) at (0,0)
    {\textbf{Map variables}\\[-0.3mm]
     $\phi=(q,\bm w)$\\
     $q$: radial grading exponent; $\bm w$: density-correction weights};
  \node[train] (map) at (4.05,0)
    {\textbf{Coordinate layer}\\[-0.3mm]
     $r_\phi,\ r_\phi'>0$\\
     integrates and normalizes positive $\rho_\phi$};
  \node[train] (solve) at (8.10,0)
    {\textbf{Galerkin equilibrium}\\[-0.3mm]
     $K(\phi)\bm d=\bm f$\\
     $K$: stiffness; $\bm f$: load; $\bm d$: spline coefficients};
  \node[train] (loss) at (12.15,0)
    {\textbf{Equilibrium energy objective}\\[-0.3mm]
     $\mathcal J_h=\tfrac12\bm d^TK\bm d$\\
     returns strain energy evaluated at equilibrium};
  \draw[flow] (parameters) -- node[above]{construct} (map);
  \draw[flow] (map) -- node[above]{assemble} (solve);
  \draw[flow] (solve) -- node[above]{evaluate} (loss);
  \draw[flow,blue!70!black] (loss.north) -- ++(0,0.75) -|
    node[pos=0.70,right]{update only $q$ and $\bm w$} (parameters.north);

  \node[eval] (frozen) at (0,-4.25)
    {\textbf{Frozen coordinate}\\[-0.3mm]
     $\phi^*$: accepted map parameters\\
     no further optimizer update};
  \node[eval] (evalsolve) at (4.05,-4.25)
    {\textbf{Evaluation solve}\\[-0.3mm]
     $K(\phi^*)\bm d^*=\bm f$\\
     same Galerkin equations on the fixed map};
  \node[eval] (metrics) at (8.10,-4.25)
    {\textbf{Reported functionals}\\[-0.3mm]
     $e_u,e_E$\\
     $K_I,K_{II}$\\
     field errors and fracture outputs};
  \node[ref] (reference) at (12.15,-4.25)
    {\textbf{Reference branch}\\[-0.3mm]
     $\bm u_{\rm ref},\bm K_{\rm ref}$\\
     used only to compute reported metrics};
  \draw[flow,teal!70!black] (map.south) -- node[midway, above, sloped]{freeze} (frozen.north);
  \draw[flow,teal!70!black] (frozen) -- node[above]{solve} (evalsolve);
  \draw[flow,teal!70!black] (evalsolve) -- node[above]{post-process} (metrics);
  \draw[flow,orange!80!black] (reference) -- node[above]{compare} (metrics);
  \node[note,blue!70!black] at (3.8,2.18) {Mechanics-only training};
  \node[note,teal!70!black] at (3.8,-2.57) {Post-training evaluation};
\end{tikzpicture}%
}
\caption{Information flow during mechanics training and post-training
evaluation.}
\label{fig:trainingevaluation}
\end{figure}
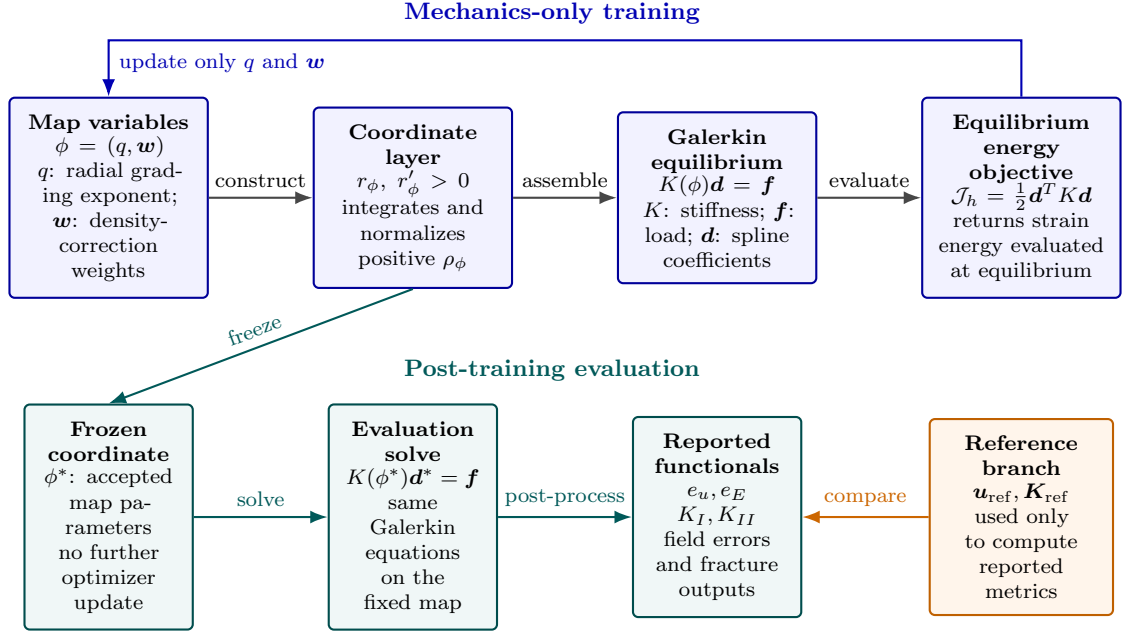

The dependency
\[
 \phi\longmapsto\rho_\phi\longmapsto r_\phi
 \longmapsto K(\phi)\longmapsto\bm d(\phi)
 \longmapsto\mathcal J_h(\phi)
\]
defines the training loop.  Figs.~\ref{fig:parameterrule} and
\ref{fig:positivedensity} realize its parameter and coordinate stages, while
Fig.~\ref{fig:mappedgalerkin} supplies the assembled equilibrium operator.
Fig.~\ref{fig:trainingevaluation} shows that the reference branch remains
disconnected until $\phi^*$ has been fixed. In that figure, $e_u$ denotes the
relative displacement error and $e_E$ the relative energy-norm error; their
precise norms are given in Section~\ref{sec:verification}.
The direct and conditional rules differ only in the source of $\phi$: the
former optimizes one parameter vector, whereas the latter evaluates
$\phi_\vartheta(\zeta)$ for a prescribed problem parameter $\zeta$.  The vector
$\bm w=(w_j)_j$ collects the density-correction weights, and $\mathcal W$
denotes their admissible set.  In either case, the network constructs the
coordinate and does not approximate the displacement or stress field.

Fig.~\ref{fig:parameterrule} distinguishes two ways of supplying the same
map parameter vector.  In the direct rule, $\phi=(q,\bm w)$ is optimized for
one mechanics problem: $q$ is the radial grading exponent, whereas
$\bm w$ modifies the positive density away from a pure power law.  In the
conditional rule, the shared network parameters $\vartheta$ convert a physical
or geometric problem parameter $\zeta$ into $\phi_\vartheta(\zeta)$, so that one rule
represents a family of singular coordinates.  The restrictions $q>0$ and
$\bm w\in\mathcal W$ are imposed before either output enters the coordinate
construction.  Both arrows therefore terminate at the same admissible
parameter interface; the subsequent Galerkin formulation is unchanged by
whether $\phi$ was identified for one problem or predicted conditionally.

In Fig.~\ref{fig:positivedensity}, $s\in[0,1]$ is the computational radius and
$\rho_\phi(s)$ is an unnormalized radial density.  The numerator accumulates
this density from the collapsed tip to $s$, while the denominator rescales the
total accumulation to unity.  Consequently,
\[
 r_\phi'(s)=\frac{\rho_\phi(s)}
 {\int_0^1\rho_\phi(\xi)\,\dd\xi}>0,
\]
which explains simultaneously the two anchors $r_\phi(0)=0$ and
$r_\phi(1)=1$ and the absence of radial folding.  Near the tip,
$\rho_\phi(s)\sim C_\phi s^{q-1}$ gives
$r_\phi(s)\sim \widetilde C_\phi s^q$; hence a physical singular term
$r^\lambda$ is pulled back to
$\widetilde C_\phi^\lambda s^{q\lambda}$.  The exponent $q$ therefore controls
the regularity gain, while $\bm w$ can redistribute resolution at finite
radius without violating positivity.  This identity is the reason that map
admissibility is built into the parameterization rather than enforced by a
penalty after training.

Fig.~\ref{fig:mappedgalerkin} shows how the one-dimensional radial map enters
the mechanics calculation.  The pullback
$\bm\chi_\phi(s,a)=r_\phi(s)(\cos a,\sin a)$ sends the entire edge $s=0$ to
the crack tip, while the two angular limits retain the two crack faces.  In
the stiffness integral, $\widehat\Omega$ is the parametric domain,
$B_\phi$ is the strain--displacement operator after this pullback, $C$ is the
constitutive tensor, and $J_\phi$ is the integration Jacobian.  Thus
$\phi$ changes $B_\phi$ and $J_\phi$, and therefore the entries of $K(\phi)$,
but it does not add basis functions or alter $C$.  Partitioning the assembled
system into free and prescribed coefficients gives the last equation in the
figure: $K_{ff}$ acts on the unknown vector $\bm d_f$, whereas
$K_{fc}\bm d_c$ transfers the imposed boundary data to the right-hand side.
This separation is important because any error reduction can be attributed to
the coordinate rather than to a larger trial space or a different material
model.

The upper branch of Fig.~\ref{fig:trainingevaluation} is the complete training
loop.  A trial $\phi$ defines $r_\phi$, which fixes the assembled matrix
$K(\phi)$; the Galerkin equation then determines the equilibrium coefficients
$\bm d(\phi)$ before
$\mathcal J_h(\phi)=\tfrac12\bm d(\phi)^TK(\phi)\bm d(\phi)$ is returned to
the optimizer.  The feedback arrow updates only $q$ and $\bm w$: the spline
coefficients are eliminated by equilibrium and are not neural-network
outputs.  After convergence, the lower branch freezes $\phi^*$, solves the
same Galerkin equations on that coordinate, and evaluates these field errors
and the SIF-extraction functionals.  The reference solution
$(\bm u_{\rm ref},\bm K_{\rm ref})$ enters only at this last comparison stage.
The disconnection of the reference branch from the training loop is essential:
it prevents exact fields or stress intensity factors from acting as hidden
labels and makes transfer to an independently evaluated mechanics problem a
genuine test of the learned coordinate.

For the scalar problem, only $q,w_1,w_2$ are trainable.  The two fixed hidden
features leave $q$ to control the radial grading exponent and use the weights to modify
the coordinate away from the tip.  Starting from the power map with zero
density weights, we first minimize \eqref{eq:objective} with respect to $q$ and
then admit bounded weight corrections.  Training is performed on a coarse
scalar discretization; the frozen map is evaluated on a finer scalar mesh and
transferred without retraining to the plane-strain problem.  Thus a change in
the vector result can be attributed to the learned coordinate rather than to a
second optimization.

For the wedge family with half-angle $\alpha$, the conditional exponent predictor uses the normalized
input \WedgeNormalizedInput{} and the bounded form in
Eq.~\eqref{eq:wedgepredictor},
\begin{align}
 x(\alpha)&=\frac{\alpha/\pi-a_x}{b_x},
 \qquad
 z_\vartheta(\alpha)=b+\sum_j v_j\tanh\!\left(\gamma_x[x(\alpha)-c_j]\right),
 \nonumber\\
 q_\vartheta(\alpha)&=q_{\min}+
 (q_{\max}-q_{\min})\,\sigma\!\left(z_\vartheta(\alpha)\right),
 \label{eq:wedgepredictor}
\end{align}
where $\sigma$ is the logistic function.  It confines every predicted grading exponent to
the positive interval $[q_{\min},q_{\max}]$.  The hidden centers and slopes are
fixed, while the bias and output weights are determined by the mean Galerkin
energy over the training wedges.  An affine rule is optimized with the same
objective, so any difference between the two predictors is caused by their
parameterization rather than by their training data.  The half-angles
\WedgeTrainAngles{} are used for training and \WedgeTestAngles{} for evaluation;
the exact exponents are withheld until both maps have been frozen.

For the nonlinear Robin family, a log-normalized Robin parameter is supplied to the
quadratic softplus predictor in Eq.~\eqref{eq:nonlinearpredictor},
\begin{align}
 z(\kappa)&=\frac{\log(\kappa/\kappa_c)}{s_\kappa},
 &q(\kappa)&=q_0+\operatorname{softplus}(a_0+a_1z+a_2z^2),
 \nonumber\\
 w_j&=w_{\max}\tanh(\widetilde w_j),
 &\rho(s;\kappa)&=s^{q(\kappa)-1}
 \exp\!\left[\sum_jw_j\tanh\!\bigl(\beta_j(s-c_j)\bigr)\right].
 \label{eq:nonlinearpredictor}
\end{align}
The logarithmic normalization resolves both ends of the parameter range, while
the softplus keeps $q(\kappa)$ above the base value $q_0$.  Three coefficients
describe this exponent rule and two bounded density weights are shared by all
parameter values.  We first determine the pure-power rule and then release the two
density weights from several independent initializations.  Because the
training and test parameter values remain unchanged, this continuation isolates the
effect of the non-power-law correction.

Let $F_h$ denote the computed radial mode and $\lambda_h$ its discrete
singular exponent, obtained from the angular eigenproblem described below.
Let $\mathcal E_c$ and $\mathcal E_e$ denote the coarse and enriched discrete
energies, $\mathcal E_0$ the identity-map energy, and $A_c$ and $A_e$ the
corresponding near-tip amplitude estimates. With fixed positive weights
$\gamma_E$, $\gamma_R$, $\gamma_A$, and $\gamma_\vartheta$, and amplitude
floor $\epsilon_A$, the reference-solution-free nonlinear objective is
\begin{align}
 \ell_m&=\log|F_h(r_m)|-\lambda_h\log r_m,
 \nonumber\\
 \delta_E&=\frac{\max(\mathcal E_c-\mathcal E_e,0)}{\mathcal E_0},
 \nonumber\\
 \mathcal L(\vartheta)&=\operatorname{mean}_{\kappa}\!\left[
 \frac{\mathcal E_c}{\mathcal E_0}+\gamma_E\delta_E
 +\gamma_R\operatorname{Var}_m(\ell_m)
 \right.
 \nonumber\\
 &\quad\left.
 +\gamma_A\left(\frac{A_c-A_e}{\max(|A_e|,\epsilon_A)}\right)^2\right]
 +\gamma_\vartheta\norm{\vartheta}^2 .
 \label{eq:nonlinearloss}
\end{align}
Here $\vartheta$ contains the exponent coefficients and, in the density-corrected model,
the density weights. The weights $\gamma_E$, $\gamma_R$, and $\gamma_A$ scale
the hierarchical-energy, near-tip radial-consistency, and amplitude terms,
respectively. The coarse--enriched discrepancy prevents the map from lowering
the coarse energy at the expense of the resolved solution,
whereas $\operatorname{Var}_m(\ell_m)$ penalizes a radial dependence that is
inconsistent with the computed eigenvalue.  The positive floor $\epsilon_A$
regularizes the amplitude ratio, and $\gamma_\vartheta\|\vartheta\|^2$ limits
unnecessary parameter growth.  All weights are fixed before optimization.
Infeasible coordinates are rejected at assembly, but neither the continuum
characteristic equation nor the reference profile enters the loss.

\subsection{Nonlinear Robin eigenfamily}

We examine the positive-density correction on a manufactured anti-plane
problem that is distinct from the linear wedge interpolation.  On the
re-entrant wedge
$\Omega_\alpha=\{(r,\theta):0<r<1,-\alpha<\theta<\alpha\}$, with
$\alpha=\NonlinearHalfAngleOverPi\pi$, the governing equation and radial
material profile are given in Eq.~\eqref{eq:nonlinearstrong}:
\begin{equation}
 -\nabla\!\cdot\!\left(\mu(r)\nabla u\right)=0,
 \qquad
 \mu(r)=\exp\!\left(c_\mu r(1-r)\right),
 \qquad c_\mu=\NonlinearModulusMultiplier{}\beta
 =\NonlinearModulusCoefficient{},\quad \beta=\NonlinearBeta{}.
 \label{eq:nonlinearstrong}
\end{equation}
Here $\mu(r)$ is the radially graded shear modulus, $\beta$ is its
dimensionless grading parameter, and $c_\mu$ fixes the strength of the radial
variation.
The two wedge faces carry equal dimensionless Robin springs $\kappa>0$.
Using the outward angular normals, the face conditions are
$-\Phi'(-\alpha)+\kappa\Phi(-\alpha)=0$ and
$\Phi'(\alpha)+\kappa\Phi(\alpha)=0$; equivalently,
$\partial_n u+\kappa u/r=0$ on both faces.
For the separated field $u(r,\theta)=F(r)\Phi(\theta)$, the angular mode is
computed from the assembled generalized eigenproblem in
Eq.~\eqref{eq:robineigen}:
\begin{equation}
 \int_{-\alpha}^{\alpha}\Phi_h'\Psi_h'\dd\theta
 +\kappa\!\sum_{\theta_f\in\{-\alpha,\alpha\}}
   \Phi_h(\theta_f)\Psi_h(\theta_f)
 =\lambda_h^2\int_{-\alpha}^{\alpha}\Phi_h\Psi_h\dd\theta .
 \label{eq:robineigen}
\end{equation}
Here $\Phi_h$ and $\Psi_h$ are the discrete angular trial and test modes, and
$\lambda_h$ is the corresponding discrete singular exponent.
Thus the angular exponent is obtained from the assembled eigenproblem during
training.  The symmetric continuum characteristic
$\lambda\tan(\lambda\alpha)=\kappa$ is evaluated only after training to report
exponent error and to construct the independent reference solution.

For $F_h(0)=0$ and $F_h(1)=1$, the mapped radial Galerkin functional in
Eq.~\eqref{eq:nonlinearradial} is assembled at every objective evaluation:
\begin{equation}
 \mathcal E_h(r_\phi;\kappa)=\frac{1}{2}\int_0^1
 \mu(r_\phi)\left[
   \frac{r_\phi}{r_\phi'}(F_{h,s})^2
   +\lambda_h^2\frac{r_\phi'}{r_\phi}F_h^2
 \right]\dd s .
 \label{eq:nonlinearradial}
\end{equation}
Here $F_h$ is the discrete radial mode and $\mathcal E_h$ is its mapped
Galerkin energy.
The training loss combines this assembled energy with hierarchical energy and
tip-amplitude defects and a near-tip equilibrium-consistency residual.
The graded modulus makes $F$ non-power-law away from the tip, so a pure-power
map and a nonzero density correction are distinguishable on this family.

For independent evaluation, let $\lambda_{\mathrm{ref}}$,
$\Phi_{\mathrm{ref}}$, and $A_{\mathrm{ref}}$ denote the continuum exponent,
angular mode, and near-tip amplitude.  We obtain the exponent by root finding
and the radial amplitude from a separately refined ODE solution.  Neither
quantity is evaluated during training.  The discrete angular projection and
the scalar amplitude error are
\begin{align}
 P_h&=\frac{\int_{-\alpha}^{\alpha}\Phi_h\Phi_{\mathrm{ref}}\dd\theta}
 {\int_{-\alpha}^{\alpha}\Phi_{\mathrm{ref}}^2\dd\theta},
 &A_h&=\frac{1}{N_p}\sum_{m=1}^{N_p}
 \frac{P_hF_h(r_m)}{r_m^{\lambda_{\mathrm{ref}}}},
 \nonumber\\
 e_A&=\frac{|A_h-A_{\mathrm{ref}}|}{|A_{\mathrm{ref}}|}.
 \label{eq:amplitudeerror}
\end{align}
Here $P_h$ is the angular projection factor, $A_h$ is the discrete near-tip
amplitude, and $r_m$, $m=1,\ldots,N_p$, are the prescribed probe radii.  In
Eq.~\eqref{eq:amplitudeerror}, $e_A$ measures the relative amplitude error with
respect to $A_{\rm ref}$ and is distinct from the mode-I and mode-II stress intensity factors
used for the elasticity problem.  Its probe radii, projection resolution, and
reference-solver settings are fixed before evaluation, and the exact profile
does not enter training.

\section{Verification problems and numerical procedures}
\label{sec:verification}

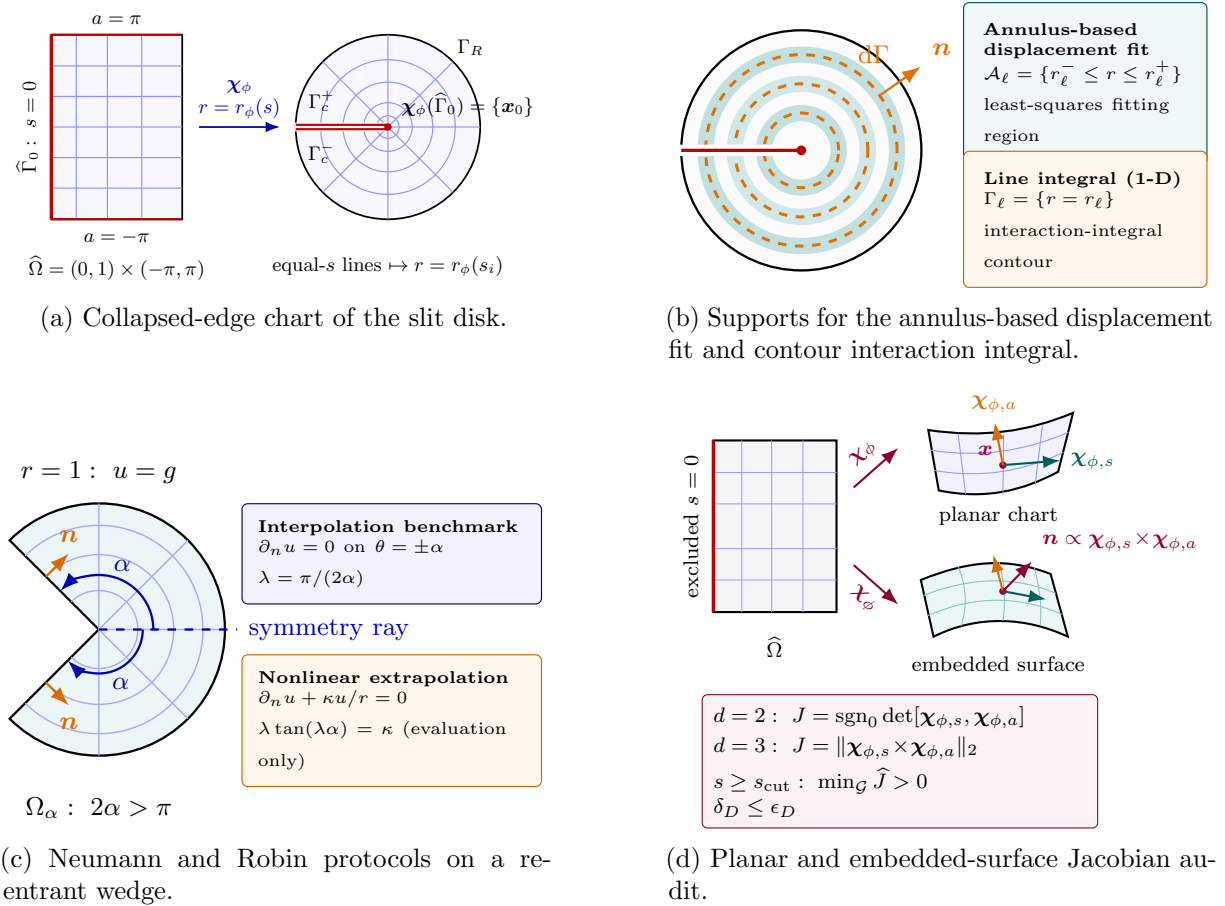
\begin{figure}[!t]
\centering
\begin{subfigure}[t]{0.45\linewidth}
\centering
\resizebox{\VerificationSchematicWidth}{!}{%
\begin{tikzpicture}[x=1cm,y=1cm,font=\footnotesize,>=Latex]
  \fill[blue!3] (-3.35,-1.48) rectangle (-1.25,1.48);
  \draw[thick] (-3.35,-1.48) rectangle (-1.25,1.48);
  \foreach \x in {-2.95,-2.45,-1.90}
    \draw[blue!42,line width=.55pt] (\x,-1.48)--(\x,1.48);
  \foreach \y in {-0.98,-0.49,0,0.49,0.98}
    \draw[blue!42,line width=.55pt] (-3.35,\y)--(-1.25,\y);
  \draw[red!75!black,line width=1.5pt] (-3.35,-1.48)--(-3.35,1.48);
  \draw[red!75!black,line width=1.0pt] (-3.35,1.48)--(-1.25,1.48);
  \draw[red!75!black,line width=1.0pt] (-3.35,-1.48)--(-1.25,-1.48);
  \node[above] at (-2.30,1.52) {$a=\pi$};
  \node[below] at (-2.30,-1.52) {$a=-\pi$};
  \node[rotate=90,above] at (-3.43,0) {$\widehat\Gamma_0:\ s=0$};
  \node[below] at (-2.30,-1.93) {$\widehat\Omega=(0,1)\times(-\pi,\pi)$};
  \draw[->,thick,blue!70!black] (-0.98,0)--(0.32,0)
    node[midway,above,align=center] {$\bm\chi_\phi$\\$r=r_\phi(s)$};

  \coordinate (C) at (2.05,0);
  \fill[blue!3] (C) circle (1.48);
  \draw[thick] (C) circle (1.48);
  \foreach \r in {0.18,0.39,0.72,1.08}
    \draw[blue!42,line width=.55pt] (C) circle (\r);
  \foreach \a in {-135,-90,-45,0,45,90,135}
    \draw[blue!42,line width=.55pt] (C)--++(\a:1.48);
  \draw[white,line width=4pt] (C)--++(180:1.50);
  \draw[red!75!black,line width=1.0pt] ($(C)+(0,0.035)$)--++(180:1.48);
  \draw[red!75!black,line width=1.0pt] ($(C)+(0,-0.035)$)--++(180:1.48);
  \fill[red!75!black] (C) circle (1.7pt);
  \node[above right,align=left] at ($(C)+(0.08,0.03)$)
    {$\bm\chi_\phi(\widehat\Gamma_0)=\{\bm x_0\}$};
  \node[above left] at ($(C)+(-0.72,0.12)$) {$\Gamma_c^+$};
  \node[below left] at ($(C)+(-0.72,-0.12)$) {$\Gamma_c^-$};
  \node[above right] at ($(C)+(0.98,1.02)$) {$\Gamma_R$};
  \node[below] at (C |- 0,-1.93) {equal-$s$ lines $\mapsto r=r_\phi(s_i)$};
\end{tikzpicture}}
\caption{Collapsed-edge chart of the slit disk.}
\label{fig:verification-crack}
\end{subfigure}\hfill
\begin{subfigure}[t]{0.45\linewidth}
\centering
\resizebox{\VerificationSchematicWidth}{!}{%
\begin{tikzpicture}[x=1cm,y=1cm,font=\footnotesize,>=Latex]
  \coordinate (C) at (-1.65,0);
  \fill[gray!3] (C) circle (1.48);
  \fill[teal!24,even odd rule] (C) circle (1.28) (C) circle (1.08);
  \fill[teal!18,even odd rule] (C) circle (0.92) (C) circle (0.72);
  \fill[teal!24,even odd rule] (C) circle (0.56) (C) circle (0.36);
  \draw[thick] (C) circle (1.48);
  \foreach \r in {0.46,0.82,1.18}
    \draw[orange!88!black,dashed,line width=1.05pt] (C) circle (\r);
  \draw[white,line width=4pt] (C)--++(180:1.50);
  \draw[red!75!black,line width=1.15pt] (C)--++(180:1.48);
  \fill[red!75!black] (C) circle (1.7pt);
  \draw[->,orange!88!black,thick] ($(C)+(35:1.18)$)--++(35:0.62)
    node[above right] {$\bm n$};
  \node[orange!88!black] at ($(C)+(53:1.48)$) {$\dd\Gamma$};
  \node[draw=teal!65!black,rounded corners=2pt,fill=teal!6,
    text width=2.5cm,align=left,inner sep=2.5mm] at (1.85,0.82)
    {\tiny{\textbf{Annulus-based displacement fit}\\
     $\mathcal A_\ell=\{r_\ell^-\le r\le r_\ell^+\}$\\
     least-squares fitting region}};
  \node[draw=orange!80!black,rounded corners=2pt,fill=orange!7,
    text width=2.5cm,align=left,inner sep=2.5mm] at (1.85,-0.85)
    {\tiny{\textbf{Line integral (1-D)}\\
     $\Gamma_\ell=\{r=r_\ell\}$\\
     interaction-integral contour}};
\end{tikzpicture}}
\caption{Supports for the annulus-based displacement fit and contour interaction integral.}
\label{fig:verification-contours}
\end{subfigure}

\medskip
\begin{subfigure}[t]{0.45\linewidth}
\centering
\resizebox{\VerificationSchematicWidth}{!}{%
\begin{tikzpicture}[x=1cm,y=1cm,font=\footnotesize,>=Latex]
  \coordinate (C) at (-1.65,0);
  \fill[teal!7] (C)--++(-135:1.55) arc (-135:135:1.55)--cycle;
  \draw[thick] (C)--++(-135:1.55) arc (-135:135:1.55)--cycle;
  \foreach \r in {0.48,0.92,1.28}
    \draw[blue!34,line width=.5pt] ($(C)+(-135:\r)$) arc (-135:135:\r);
  \foreach \a in {-90,-45,0,45,90}
    \draw[blue!34,line width=.5pt] (C)--++(\a:1.55);
  \draw[blue!75!black,dashed,thick] (C)--++(0:1.70)
    node[right] {symmetry ray};
  \draw[->,blue!70!black,thick] ($(C)+(0:0.67)$) arc (0:135:0.67);
  \draw[->,blue!70!black,thick] ($(C)+(0:0.54)$) arc (0:-135:0.54);
  \node[blue!70!black] at ($(C)+(70:0.84)$) {$\alpha$};
  \node[blue!70!black] at ($(C)+(-68:0.72)$) {$\alpha$};
  \draw[->,orange!85!black,thick] ($(C)+(135:0.92)$)--++(45:0.40)
    node[above] {$\bm n$};
  \draw[->,orange!85!black,thick] ($(C)+(-135:0.92)$)--++(-45:0.40)
    node[below] {$\bm n$};
  \node[inner sep=1pt] at (-1.65,1.93) {$r=1:\ u=g$};
  \node[below] at (-1.65,-1.93) {$\Omega_\alpha:\ 2\alpha>\pi$};
  \node[draw=blue!60!black,rounded corners=2pt,fill=blue!5,
    text width=3.28cm,align=left,inner sep=2mm] at (1.95,0.93)
    {\tiny{\textbf{Interpolation benchmark}\\
     $\partial_nu=0$ on $\theta=\pm\alpha$\\
     $\lambda=\pi/(2\alpha)$}};
  \node[draw=orange!80!black,rounded corners=2pt,fill=orange!7,
    text width=3.28cm,align=left,inner sep=2mm] at (1.95,-1.12)
    {\tiny{\textbf{Nonlinear extrapolation}\\
     $\partial_nu+\kappa u/r=0$\\
     $\lambda\tan(\lambda\alpha)=\kappa$ (evaluation only)}};
\end{tikzpicture}}
\caption{Neumann and Robin protocols on a re-entrant wedge.}
\label{fig:verification-wedge}
\end{subfigure}\hfill
\begin{subfigure}[t]{0.45\linewidth}
\centering
\resizebox{\VerificationSchematicWidth}{!}{%
\begin{tikzpicture}[x=1cm,y=1cm,font=\scriptsize,>=Latex]
  \fill[gray!8] (-3.50,-1.15) rectangle (-1.85,1.15);
  \draw[thick] (-3.50,-1.15) rectangle (-1.85,1.15);
  \foreach \x in {-3.12,-2.72,-2.30}
    \draw[blue!38] (\x,-1.15)--(\x,1.15);
  \foreach \y in {-0.72,-0.30,0.30,0.72}
    \draw[blue!38] (-3.50,\y)--(-1.85,\y);
  \draw[red!75!black,line width=1.3pt] (-3.50,-1.15)--(-3.50,1.15);
  \node[rotate=90,above] at (-3.57,0) {excluded $s=0$};
  \node[below] at (-2.68,-1.34) {$\widehat\Omega$};

  \draw[->,purple!75!black,thick] (-1.62,0.52)--(-1.02,1.05)
    node[midway,above,sloped] {$\bm\chi_\phi$};
  \draw[->,purple!75!black,thick] (-1.62,-0.52)--(-1.02,-1.05)
    node[midway,below,sloped] {$\bm\chi_\phi$};

  \fill[blue!5] (-0.48,0.48) .. controls (0.15,0.35) and (0.65,0.48) .. (1.12,0.73)
    -- (1.32,1.52) .. controls (0.72,1.27) and (0.12,1.18) .. (-0.62,1.29) -- cycle;
  \draw[thick] (-0.48,0.48) .. controls (0.15,0.35) and (0.65,0.48) .. (1.12,0.73)
    -- (1.32,1.52) .. controls (0.72,1.27) and (0.12,1.18) .. (-0.62,1.29) -- cycle;
  \draw[blue!42] (-0.55,0.75) .. controls (0.12,0.63) and (0.68,0.75) .. (1.19,0.99);
  \draw[blue!42] (-0.59,1.02) .. controls (0.10,0.91) and (0.70,1.02) .. (1.26,1.26);
  \draw[blue!42] (-0.12,0.43)--(-0.22,1.24);
  \draw[blue!42] (0.38,0.42)--(0.39,1.23);
  \draw[blue!42] (0.81,0.55)--(0.91,1.35);
  \fill[purple!80!black] (0.38,0.82) circle (1.4pt) node[above left] {$\bm x$};
  \draw[->,teal!75!black,thick] (0.38,0.82)--(1.15,0.88)
    node[right] {$\bm\chi_{\phi,s}$};
  \draw[->,orange!85!black,thick] (0.38,0.82)--(0.27,1.39)
    node[above] {$\bm\chi_{\phi,a}$};
  \node[below] at (0.32,0.4) {planar chart};

  \fill[teal!7] (-0.55,-1.48) .. controls (0.05,-1.18) and (0.70,-1.22) .. (1.25,-1.48)
    -- (1.02,-0.64) .. controls (0.45,-0.42) and (-0.05,-0.45) .. (-0.72,-0.72) -- cycle;
  \draw[thick] (-0.55,-1.48) .. controls (0.05,-1.18) and (0.70,-1.22) .. (1.25,-1.48)
    -- (1.02,-0.64) .. controls (0.45,-0.42) and (-0.05,-0.45) .. (-0.72,-0.72) -- cycle;
  \draw[teal!48] (-0.63,-1.20) .. controls (0.02,-0.91) and (0.63,-0.94) .. (1.17,-1.18);
  \draw[teal!48] (-0.69,-0.93) .. controls (0.00,-0.66) and (0.58,-0.68) .. (1.09,-0.90);
  \draw[teal!48] (-0.18,-1.34) .. controls (-0.08,-1.02) and (-0.13,-0.75) .. (-0.22,-0.54);
  \draw[teal!48] (0.38,-1.27) .. controls (0.43,-0.94) and (0.38,-0.67) .. (0.33,-0.48);
  \draw[teal!48] (0.88,-1.33) .. controls (0.91,-1.00) and (0.85,-0.76) .. (0.80,-0.54);
  \fill[purple!80!black] (0.38,-0.87) circle (1.4pt);
  \draw[->,teal!75!black,thick] (0.38,-0.87)--(0.96,-0.98);
  \draw[->,orange!85!black,thick] (0.38,-0.87)--(0.25,-0.38);
  \draw[->,purple!75!black,thick] (0.38,-0.87)--(0.77,-0.47)
    node[above right] {$\bm n\propto\bm\chi_{\phi,s}\!\times\!\bm\chi_{\phi,a}$};
  \node[below] at (0.32,-1.61) {embedded surface};

  \node[draw=purple!65!black,rounded corners=2pt,fill=purple!5,
    text width=5cm,align=left,inner sep=1.5mm] at (-1,-3.15)
    {$d=2:\ J=\operatorname{sgn}_0\det[\bm\chi_{\phi,s},\bm\chi_{\phi,a}]$\\[0.8mm]
     $d=3:\ J=\|\bm\chi_{\phi,s}\!\times\!\bm\chi_{\phi,a}\|_2$\\[0.8mm]
     $s\ge s_{\rm cut}:\ \min_{\mathcal G}\widehat J>0$\\
     $\delta_D\le\epsilon_D$};
\end{tikzpicture}}
\caption{Planar and embedded-surface Jacobian audit.}
\label{fig:verification-chart}
\end{subfigure}
\caption{Verification geometries, supports, boundary conditions, and chart
diagnostics.}
\label{fig:verificationproblems}
\end{figure}

The numerical tests are arranged so that the coordinate is the only changing
component whenever an accuracy gain is attributed to the map.  The same radial
chart is therefore used in matched scalar and elastic operators, and the
resulting elastic field is passed to an annulus-based displacement fit and an independent contour
interaction integral for $(K_I,K_{II})$.  Refinement and parameter transfer then show
whether the gain persists beyond one discretization, while Jacobian and
computational cost at a prescribed error tolerance determine whether it is obtained by an admissible
and computationally meaningful map.  The geometries and evaluation supports
are shown in
Fig.~\ref{fig:verificationproblems}; the wedge protocols and chart audit are
detailed in Figs.~\ref{fig:verification-wedge} and
\ref{fig:verification-chart}, respectively.  The formulations below keep the
physical operator, trial-space size, and output procedure fixed whenever the
effect of the coordinate is compared. In the crack panel, $\Gamma_c^+$ and
$\Gamma_c^-$ label the upper and lower crack faces and $\Gamma_R$ the outer
boundary. In the chart-audit panel, $\mathcal G$
is the finite audit grid outside the collapsed-edge cutoff $s_{\rm cut}$,
$\widehat J=J/\max_{\mathcal G}J$ is the normalized Jacobian, and $\delta_D$
is the relative analytic--finite-difference derivative discrepancy tested
against the prescribed tolerance $\epsilon_D$.

\subsection{Matched scalar and elastic Galerkin formulations}

To compare the coordinates on the same physical problem and discrete space,
we define the slit disk and radial chart by
\[
 \begin{aligned}
 \Omega_R&=\{(r,a):0<r<R,\ -\pi<a<\pi\},&
 \widehat\Omega&=(0,1)\times(-\pi,\pi),\\
 \bm\chi_\phi(s,a)&=r_\phi(s)(\cos a,\sin a).&&
 \end{aligned}
\]
The two limits $a=\pm\pi$ represent distinct crack faces, although
$\bm\chi_\phi(0,a)$ collapses to one physical tip.  For a scalar field,
\[
 \nabla_\phi v=
 \frac{v_{,s}}{r_\phi'}\bm e_r+
 \frac{v_{,a}}{r_\phi}\bm e_a,\qquad
 J_\phi=r_\phi r_\phi'.
\]
Here $\bm e_r$ and $\bm e_a$ are the radial and angular unit vectors,
respectively.
Consequently, imposing the trace from Eq.~\eqref{eq:scalartraining} gives the
single variational statement
\begin{equation}
 \begin{aligned}
 &\text{find }u_h-g_h\in V_h^0(r_\phi)
 \text{ such that}\quad
 A_\phi^{\rm sc}(u_h,v_h)=0
 &&\forall v_h\in V_h^0(r_\phi),\\
 &A_\phi^{\rm sc}(u,v)=
 \int_{-\pi}^{\pi}\!\!\int_0^1
 \left(
 \frac{r_\phi}{r_\phi'}u_{,s}v_{,s}
 +\frac{r_\phi'}{r_\phi}u_{,a}v_{,a}
 \right)\dd s\,\dd a .
 \end{aligned}
 \label{eq:scalarweakform}
\end{equation}
The coordinate enters Eq.~\eqref{eq:scalarweakform} through the two reciprocal
metric factors.  The physical differential operator and the spline basis are
fixed for all coordinate choices.

For plane strain, define
$\bm\varepsilon_\phi(\bm v)=\operatorname{sym}(\nabla_\phi\bm v)$ and
$\bm\sigma_\phi(\bm v)=\bm C(E,\nu):\bm\varepsilon_\phi(\bm v)$, where
$\bm C(E,\nu)$ is the isotropic plane-strain elasticity tensor determined by
Young's modulus $E$ and Poisson's ratio $\nu$.  The physical benchmark in
Fig.~\ref{fig:verification-crack} is fully specified by
\begin{equation}
 \left.
 \begin{aligned}
 -\nabla\!\cdot\!\bm\sigma(\bm u)&=\bm0 &&\text{in }\Omega_R,\\
 \bm\sigma(\bm u)\bm n&=\bm0 &&\text{on }\Gamma_c^+\cup\Gamma_c^-,\\
 \bm u&=\bm g(K_I,K_{II},T)
 &&\text{on }\Gamma_R,\\
 (R,E,\nu,K_I,K_{II},T)&=
 (\VectorBenchmarkRadius{},\VectorYoungModulus{},\VectorPoissonRatio{},
 \VectorTargetKI{},\VectorTargetKII{},\VectorTargetTStress{}).
 \end{aligned}
 \right\}
 \label{eq:vectorbenchmark}
\end{equation}
Here $R$ is the disk radius, $E$ and $\nu$ are Young's modulus and Poisson's
ratio, $K_I$ and $K_{II}$ are the mode-I and mode-II stress intensity factors,
and $T$ is the nonsingular $T$-stress.  Moreover, $\Gamma_c^\pm$ are the two
crack faces, $\Gamma_R=\{r=R\}$, and
$\bm g=\bm u_{K_I,K_{II}}+\bm u_T$ is the mixed-mode Williams trace plus its
bounded nonsingular $T$-stress term.  The Galerkin counterpart of
Eq.~\eqref{eq:vectorbenchmark} is
\begin{equation}
 \begin{aligned}
 &\text{find }\bm u_h-\bm g_h\in[ V_h^0(r_\phi)]^2
 \text{ such that}\quad
 A_\phi^{\rm el}(\bm u_h,\bm v_h)=0
 &&\forall\bm v_h\in[ V_h^0(r_\phi)]^2,\\
 &A_\phi^{\rm el}(\bm u,\bm v)=
 \int_{\widehat\Omega}
 \bm\varepsilon_\phi(\bm v):\bm C:
 \bm\varepsilon_\phi(\bm u)\,J_\phi\dd s\,\dd a .
 \end{aligned}
 \label{eq:elasticweakform}
\end{equation}
Equation~\eqref{eq:elasticweakform} yields the assembled system in
Eq.~\eqref{eq:discrete}.  All baselines consequently use the same elasticity
law and coefficient count; their strain--geometry factors differ through
$r_\phi$ and $r_\phi'$.

The matched-size comparison is the four-member set
\[
 \mathfrak B=\left\{
 \begin{array}{lll}
 \text{identity map:}   &r(s)=s,   &\xi_i=i/n,\\
 \text{graded:}  &r(s)=s,   &\xi_i=(i/n)^2,\\
 \text{fixed:}   &r(s)=s^2, &\xi_i=i/n,\\
 \text{learned:} &r(s)=r_{\phi^*}(s),&\xi_i=i/n
 \end{array}\right\},\qquad i=1,\ldots,n-1.
\]
Thus degree, angular and radial span counts, quadrature, boundary projection,
and coefficient count are invariant over $\mathfrak B$.  The radial open knot
vector and the three field metrics are
\[
 \Xi_r^{(g)}=[0^{p+1},\{(i/n)^g\}_{i=1}^{n-1},1^{p+1}],
 \qquad
 e_u=\frac{\norm{\bm u_h-\bm u_{\rm ref}}_{L^2(\Omega)}}
 {\norm{\bm u_{\rm ref}}_{L^2(\Omega)}},
\]
\[
 e_\sigma=\frac{\norm{\bm\sigma_h-\bm\sigma_{\rm ref}}_{L^2(\Omega)}}
 {\norm{\bm\sigma_{\rm ref}}_{L^2(\Omega)}},\qquad
 e_E=\left(
 \frac{A^{\rm el}(\bm u_h-\bm u_{\rm ref},\bm u_h-\bm u_{\rm ref})}
 {A^{\rm el}(\bm u_{\rm ref},\bm u_{\rm ref})}
 \right)^{1/2},\qquad
 \kappa_2=\frac{\lambda_{\max}(K_{ff})}{\lambda_{\min}(K_{ff})}.
\]
Thus $e_u$, $e_\sigma$, and $e_E$ are the relative displacement, stress, and
energy-norm errors, respectively, while $\kappa_2$ is the spectral condition
number of the free--free stiffness block. Here $p$ is the spline degree and
$g$ is the knot-grading exponent. Independent physical quadrature evaluates
all norms. The scalar radial
problems use degree-$p$ $C^0$ spaces, whereas elasticity uses the
tensor-product B-spline space in Eq.~\eqref{eq:trialspace}.

\subsection{Evaluation of stress intensity factors}

To evaluate the stress intensity factors, let $\bm U_I$ and $\bm U_{II}$
denote Williams displacement fields normalized to unit SIF,
$\{\bm P_m\}_{m=1}^{M}$ a bounded polynomial basis, and
$\mathcal A_\ell=\{r_\ell^-\leq r\leq r_\ell^+\}$ the physical fitting
annuli in Fig.~\ref{fig:verification-contours}.  The SIFs are extracted by
the following weighted least-squares fit of the Williams basis over each annulus
\begin{equation}
 \begin{aligned}
 (K_{I,\ell}^{h},K_{II,\ell}^{h},\bm c_\ell)
 &=\underset{k_I,k_{II},\bm c}{\operatorname{argmin}}\,
 \mathcal R_\ell(k_I,k_{II},\bm c),\\
 \mathcal R_\ell
 &=\int_{\mathcal A_\ell}
 \left\|
 \bm u_h-k_I\bm U_I-k_{II}\bm U_{II}
 -\sum_{m=1}^{M}c_m\bm P_m
 \right\|_2^2\dd\Omega,\\
 \ell&=1,\ldots,N_{\mathcal A},\qquad
 N_{\mathcal A}=\InteractionContourCount{}.
 \end{aligned}
 \label{eq:annularfit}
\end{equation}
Equation~\eqref{eq:annularfit} gives one mixed-mode estimate for each annulus.
Its dimensionless error is
\begin{equation}
 e_{K,\ell}=\frac{
 \left\|(K_{I,\ell}^{h},K_{II,\ell}^{h})-(K_I,K_{II})\right\|_2}
 {\left\|(K_I,K_{II})\right\|_2},\qquad
 e_K^{\max}=\max_{1\leq\ell\leq N_{\mathcal A}}e_{K,\ell}.
 \label{eq:siferror}
\end{equation}
In Eq.~\eqref{eq:siferror}, variation of $e_{K,\ell}$ with $\ell$ measures sensitivity to the fitting
annulus, whereas $e_K^{\max}$ is the maximum over all prescribed annuli.

A second functional provides an independent evaluation without the polynomial terms used in the annular fit.  For the computed state $(1)$ and a unit auxiliary Williams
state $(2)$, the interaction integral \cite{Yau1980MixedMode} is
\begin{equation}
 \mathcal I^{(1,2)}=\int_\Gamma\left[
 W^{(1,2)}n_x
 -t_i^{(1)}u_{i,x}^{(2)}
 -t_i^{(2)}u_{i,x}^{(1)}
 \right]\dd\Gamma,
 \label{eq:interaction}
\end{equation}
where $W^{(1,2)}=\sigma_{ij}^{(1)}\varepsilon_{ij}^{(2)}$,
$t_i=\sigma_{ij}n_j$, and $E'=E/(1-\nu^2)$ is the plane-strain effective
modulus.  For homogeneous isotropic plane
strain, Eqs.~\eqref{eq:interaction} and \eqref{eq:interactionrelation} are
linked by
\begin{equation}
 \mathcal I^{(1,2)}=\frac{2}{E'}
 \left(K_I^{(1)}K_I^{(2)}+K_{II}^{(1)}K_{II}^{(2)}\right).
 \label{eq:interactionrelation}
\end{equation}
Unit auxiliary modes give $K_{m,\ell}^{\mathcal I}=E'\mathcal I_{m,\ell}/2$.
Writing $\bm K_\ell^{\mathcal I}=(K_{I,\ell}^{\mathcal I},
K_{II,\ell}^{\mathcal I})$, contour invariance and target error are measured
independently as
\[
 \delta_{\rm path}=
 \frac{\max_\ell\left\|\bm K_\ell^{\mathcal I}
 -\overline{\bm K}^{\mathcal I}\right\|_2}
 {\left\|\overline{\bm K}^{\mathcal I}\right\|_2},\qquad
 e_{K,\mathcal I}^{\max}=
 \max_\ell\frac{\left\|\bm K_\ell^{\mathcal I}
 -(\InteractionTargetKI{},\InteractionTargetKII{})\right\|_2}
 {\left\|(\InteractionTargetKI{},\InteractionTargetKII{})\right\|_2},
 \qquad
 \overline{\bm K}^{\mathcal I}=\frac{1}{N_\Gamma}
 \sum_{\ell=1}^{N_\Gamma}\bm K_\ell^{\mathcal I}.
\]
The contour functional uses only $\nabla\bm u_h$, $(E,\nu)$, and the auxiliary
state.  It is restricted to the straight homogeneous crack considered here;
curved cracks and bimaterial interfaces require the corresponding domain
integrals \cite{Chiaramonte2015Curved}.

\subsection{Refinement, transfer, and admissibility tests}

To distinguish regularity transfer from an isolated-grid effect, the
refinement slope of a metric $e$ at fixed degree $p$ is evaluated on the
selected finest index set $\mathcal F$,
\[
 \widehat\beta_e=\underset{\beta,c}{\operatorname{argmin}}
 \sum_{i\in\mathcal F}
 \left[\log e_i-c-\beta\log(\mathrm{DOF}_i)\right]^2.
\]
An independent $(p,C^{p-1})$ sweep then tests whether the observed advantage
persists under changes of the spline space, without reusing this fitted slope.

The wedge interpolation benchmark has an affine angle-to-grading-exponent relation.  Its
disjoint mechanics-only training and test problems are
\begin{equation}
 \begin{aligned}
 \mathcal A_{\rm tr}/\pi&=\WedgeTrainAngles{},&
 \mathcal A_{\rm te}/\pi&=\WedgeTestAngles{},\\
 \mathcal P_\vartheta&=\WedgeParameterBounds{}^{\WedgeTrainableCount{}},&
 \mathcal J_{\rm tr}(\vartheta)
 &=\frac{1}{|\mathcal A_{\rm tr}|}
 \sum_{\alpha\in\mathcal A_{\rm tr}}
 \mathcal J_h(q_\vartheta(\alpha);\alpha)
 +\WedgeRegularization{}\norm{\vartheta}_2^2,\\
 \vartheta^*&=\underset{\vartheta\in\mathcal P_\vartheta}
 {\operatorname{argmin}}\,\mathcal J_{\rm tr}(\vartheta),&
 (a_0^*,a_1^*)&=\underset{a_0,a_1}{\operatorname{argmin}}\,
 \frac{\sum_{\alpha\in\mathcal A_{\rm tr}}
 \mathcal J_h(a_0+a_1\alpha/\pi;\alpha)}
 {|\mathcal A_{\rm tr}|}.
 \end{aligned}
 \label{eq:wedgeprotocol}
\end{equation}
All $\mathcal A_{\rm te}$ lie inside
$[\min\mathcal A_{\rm tr},\max\mathcal A_{\rm tr}]$; hence
Eq.~\eqref{eq:wedgeprotocol} tests interpolation and optimization stability.
The affine mechanics model is retained because the manufactured
$\alpha\mapsto q$ relation is affine.

The nonlinear Robin split, in contrast, imposes two-sided extrapolation:
\begin{equation}
 \begin{aligned}
 \mathcal K_{\rm tr}&=\NonlinearTrainingControls{},&
 \mathcal K_{\rm te}&=\NonlinearTestControls{},\\
 \NonlinearTestMinimum{}&<\min\mathcal K_{\rm tr}
 \leq\kappa\leq\max\mathcal K_{\rm tr}<\NonlinearTestMaximum{},
 &\kappa&\in\mathcal K_{\rm tr},\\
 \mathcal L_{\rm tr}(\vartheta)
 &=\frac{1}{|\mathcal K_{\rm tr}|}
 \sum_{\kappa\in\mathcal K_{\rm tr}}\mathcal L_\kappa(\vartheta),&
 \vartheta^*&=\underset{\vartheta}{\operatorname{argmin}}\,
 \mathcal L_{\rm tr}(\vartheta).
 \end{aligned}
 \label{eq:nonlinearsplit}
\end{equation}
The pure-power and density-corrected models share $\mathcal K_{\rm tr}$,
$\mathcal K_{\rm te}$, and the exponent predictor in
Eq.~\eqref{eq:nonlinearsplit}.  Their difference therefore measures the
contribution of the non-power-law radial-density correction.  Separate perturbations of
the radial-map normalization and Galerkin quadrature orders determine whether that
difference is numerical.

For the finite quadrature/sample set $\mathcal Q$ and implemented chart
$\bm\chi_\phi$, let $\bm x\in\mathcal Q$ denote a sample point and
$s(\bm x)$ its radial coordinate.  We write $D_{\rm FD}\bm\chi_\phi$ for the
finite-difference derivative, $J_{\rm ref}>0$ for a chart-dependent Jacobian
normalization, and $\epsilon_{\rm FD}>0$ for the denominator floor.  The
acceptance vector is
\begin{equation}
 \begin{aligned}
 \bm G(\phi)&=(G_\rho,G_b,G_J,G_D),\\
 G_\rho&=\min_{\bm x\in\mathcal Q}\rho_\phi(s(\bm x)),&
 G_b&=\max_{\zeta}
 \{|r_\phi(0;\zeta)|+|r_\phi(1;\zeta)-1|\},\\
 G_J&=\min_{\bm x\in\mathcal Q}
 \frac{J_\phi(\bm x)}{J_{\rm ref}(\bm x)},&
 G_D&=\max_{\bm x\in\mathcal Q}
 \frac{\|D\bm\chi_\phi(\bm x)-D_{\rm FD}\bm\chi_\phi(\bm x)\|_F}
 {\max(\|D_{\rm FD}\bm\chi_\phi(\bm x)\|_F,\epsilon_{\rm FD})}.
 \end{aligned}
 \label{eq:admissibilityvector}
\end{equation}
Here $G_\rho$, $G_b$, $G_J$, and $G_D$ measure, respectively, the minimum
density, boundary-anchor defect, normalized Jacobian lower bound, and relative
derivative discrepancy.  For a planar chart,
$J_\phi=\operatorname{sgn}_0\det[\bm\chi_{\phi,s},\bm\chi_{\phi,a}]$,
where $\operatorname{sgn}_0$ fixes the reference orientation; for an embedded
surface, $J_\phi=\|\bm\chi_{\phi,s}\times\bm\chi_{\phi,a}\|_2$.
The positive $J_{\rm ref}$ supplies only chart-specific normalization and is
$\max_{\mathcal G}J_\phi$ in the implemented-chart audit below. The prescribed
anchor and derivative tolerances are $\epsilon_b$ and $\epsilon_D$; hence the
required conditions are $G_\rho>0$, $G_b\leq\epsilon_b$, $G_J>0$, and
$G_D\leq\epsilon_D$.  Equation~\eqref{eq:admissibilityvector} provides a
finite-sample admissibility test away from the collapsed edge.  Global
injectivity is not implied by this test.  Let $\mathcal I_{\rm tip}$ be the indices of the basis functions excluded at
the collapsed edge. The corresponding admissible trial space $V_h^{\rm adm}$
is assembled directly as
\[
 V_h^{\rm adm}=\operatorname{span}\{N_A:A\notin\mathcal I_{\rm tip}\},
 \qquad
 (K_h)_{AB}=A_\phi(N_B,N_A),\quad A,B\notin\mathcal I_{\rm tip},
\]
so basis functions excluded at the collapsed edge never enter the assembled matrix.

For the cost comparison, let $\mathcal C_M$ contain the computed
configurations of method $M$ and let $\tau_E$ be the common energy-error upper
bound.  The computational cost at this tolerance and its amortization are defined by
\begin{equation}
 \begin{aligned}
 c_M^*&=\underset{c\in\mathcal C_M}{\operatorname{argmin}}
 \left\{\operatorname{DOF}(c):e_E(c)\leq\tau_E\right\},\\
 T_M(N)&=T_M^{\rm setup}+N\,T_M^{\rm online},\\
 N_{\rm BE}^{L/A}&=
 \left\lceil
 \frac{T_L^{\rm setup}-T_A^{\rm setup}}
 {T_A^{\rm online}-T_L^{\rm online}}
 \right\rceil_{+}.
 \end{aligned}
 \label{eq:costprotocol}
\end{equation}
Equation~\eqref{eq:costprotocol} selects the smallest tested discretization, in terms of DOF, satisfying the prescribed error tolerance.  Separating
$T^{\rm setup}$ from $T^{\rm online}$ distinguishes the cost of constructing a
map from that of reusing it.  The quantity $N_{\rm BE}^{L/A}$ is the
learned-versus-adaptive break-even solve count, with $L$ and $A$ identifying
the two methods.  Online time includes assembly and
a CG solution of the symmetrically diagonally scaled system; the adaptive reference uses the same accounting.

\section{Numerical results and discussion}
\label{sec:results}

\subsection{Straight-crack problem}

\begin{figure}[htp]
\centering
\begin{subfigure}[t]{0.32\linewidth}\centering
\includegraphics[width=\linewidth]{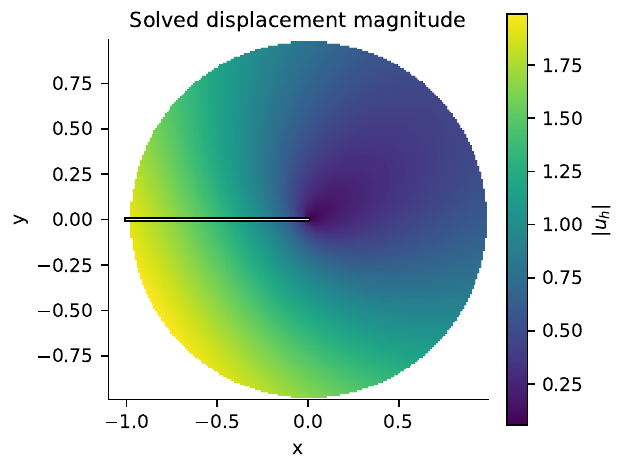}
\caption{Displacement magnitude.}
\end{subfigure}\hfill
\begin{subfigure}[t]{0.32\linewidth}\centering
\includegraphics[width=\linewidth]{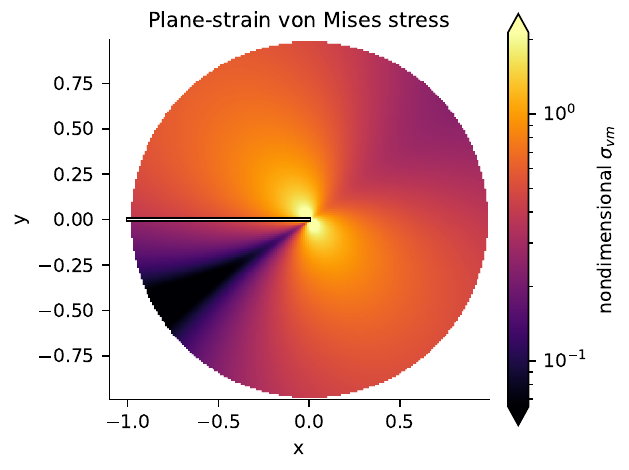}
\caption{von Mises stress.}
\end{subfigure}\hfill
\begin{subfigure}[t]{0.32\linewidth}\centering
\includegraphics[width=\linewidth]{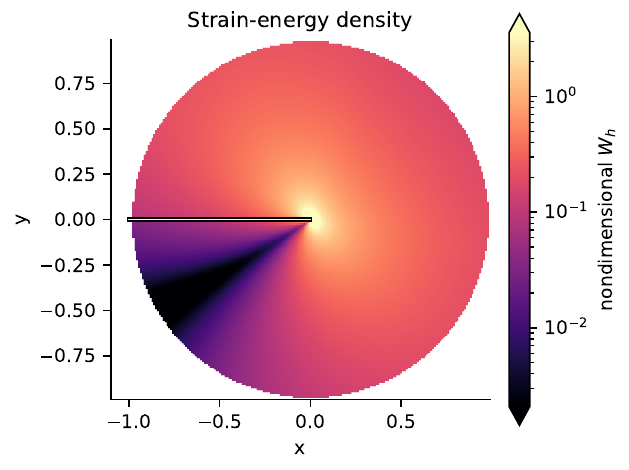}
\caption{Strain-energy density.}
\end{subfigure}

\medskip
\begin{subfigure}[t]{0.32\linewidth}\centering
\includegraphics[width=\linewidth]{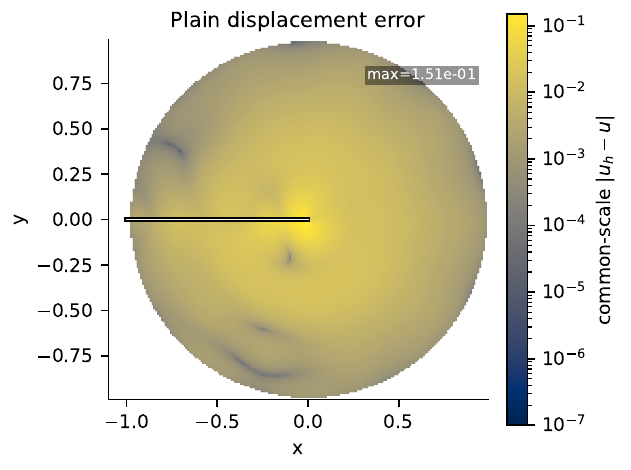}
\caption{Identity-map displacement error.}
\end{subfigure}\hfill
\begin{subfigure}[t]{0.32\linewidth}\centering
\includegraphics[width=\linewidth]{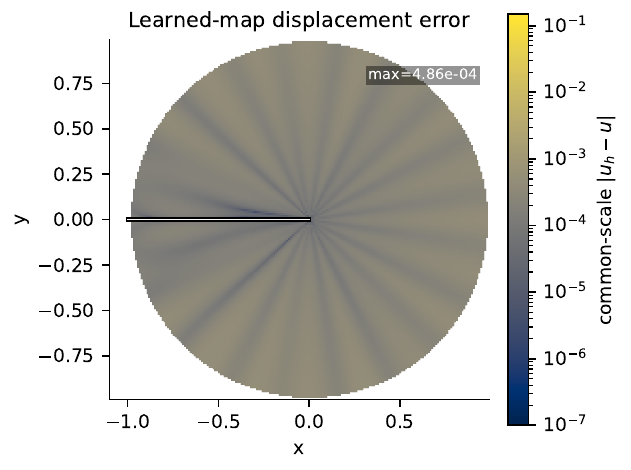}
\caption{Learned-map displacement error.}
\end{subfigure}\hfill
\begin{subfigure}[t]{0.32\linewidth}\centering
\includegraphics[width=\linewidth]{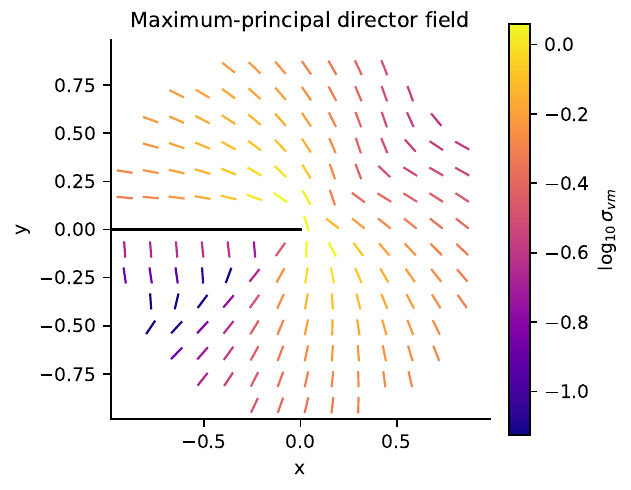}
\caption{Principal-stress directions.}
\end{subfigure}
\caption{Plane-strain displacement, stress, and error fields.}
\label{fig:vectorfields}
\end{figure}

\begin{figure}[htp]
\centering
\begin{subfigure}[t]{0.45\linewidth}\centering
\includegraphics[width=\linewidth]{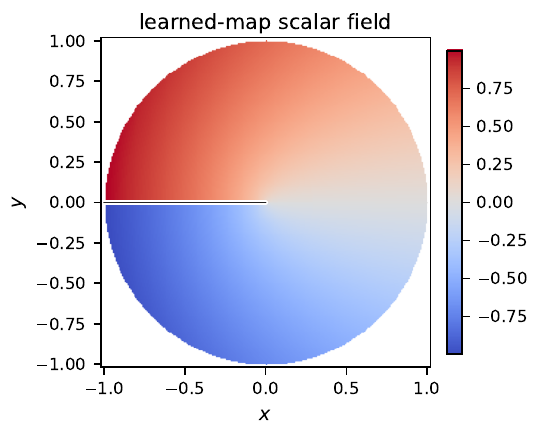}
\caption{Mapped Galerkin field.}
\end{subfigure}\hfill
\begin{subfigure}[t]{0.45\linewidth}\centering
\includegraphics[width=\linewidth]{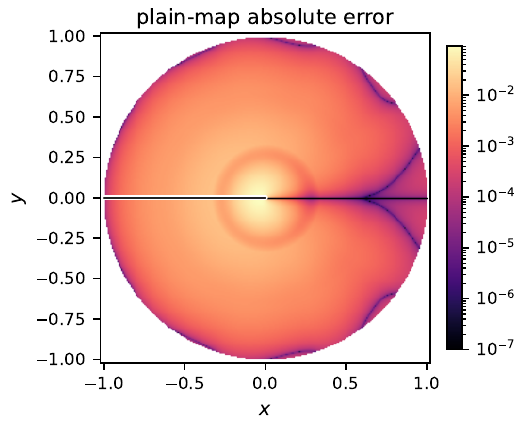}
\caption{Identity-map absolute error.}
\end{subfigure}

\medskip
\begin{subfigure}[t]{0.45\linewidth}\centering
\includegraphics[width=\linewidth]{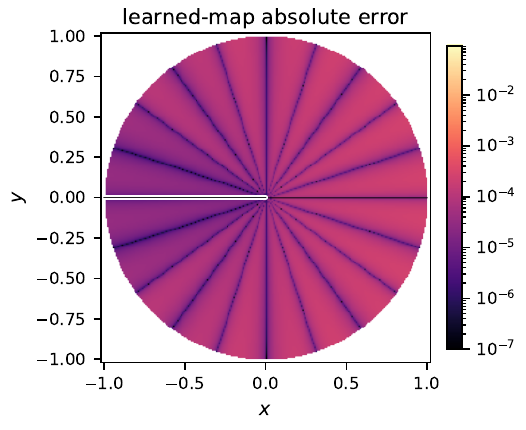}
\caption{Learned-map absolute error.}
\end{subfigure}\hfill
\begin{subfigure}[t]{0.45\linewidth}\centering
\includegraphics[width=\linewidth]{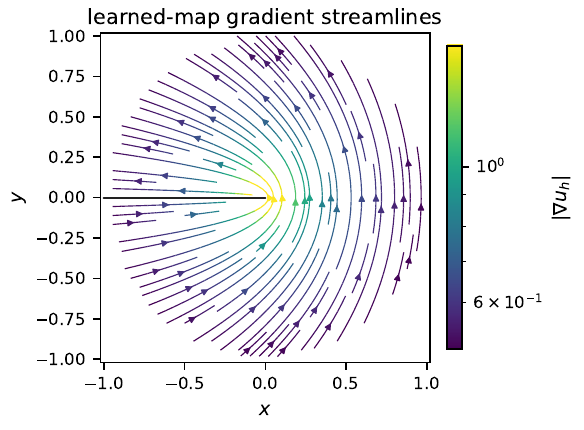}
\caption{Gradient streamlines.}
\end{subfigure}
\caption{Scalar solution, absolute-error fields, and gradient streamlines.}
\label{fig:scalarfields}
\end{figure}

At equal spline dimensions, the learned coordinate reduces the field errors
because the pullback $u(r_\phi(s),a)\sim s^{q\lambda}$ changes the resolved
radial regularity, not the size of the approximation space.  Its maximum
annulus-based displacement-fit SIF error is
\VectorLearnedSIFCorrelationError{}, and the corresponding improvement factor
over the identity map is \VectorSIFCorrelationGain{}.  The minimum eigenvalue of $K_{ff}$ is positive
(\LearnedMinimumEigenvalue{}), which confirms positive definiteness of the reduced system after imposing the essential boundary conditions.  The sampled Jacobian minimum is
\LearnedMinimumSampledJacobian{}; its small value results from the prescribed
radial collapse at the tip.  The normalized Jacobian tests for the complete
charts are considered below.  Table~\ref{tab:equal} reports these matched-size
results; its SIF error is the maximum over the displacement-fitting annuli,
and $\kappa_2(K_{ff})$ is the unpreconditioned condition number.  The vector and
scalar fields are shown in Figs.~\ref{fig:vectorfields} and
\ref{fig:scalarfields}; the identity-map and learned-map scalar errors use a common
color scale.

\begin{table}[t]
\centering
\ResultTableFont
\caption{Equal-DOF plane-strain Galerkin results.}
\label{tab:equal}
\begin{adjustbox}{max width=\linewidth}
\begin{tabular}{lrrrrrr}
\toprule
Radial mapping & DOF & rel. $L^2$ & rel. energy norm & rel. stress $L^2$ & max. annulus-based $e_K$ & $\kappa(K_{ff})$ \\
\midrule
Identity map & $156$ & $1.192\times 10^{-2}$ & $0.1684$ & $0.137$ & $8.589\times 10^{-1}$ & $100.5$ \\
Radially graded knots & $156$ & $4.569\times 10^{-3}$ & $0.09972$ & $0.08148$ & $7.7\times 10^{-1}$ & $175.3$ \\
Prescribed power map $q=2$ & $156$ & $2.425\times 10^{-4}$ & $0.00504$ & $0.004911$ & $1.015\times 10^{-4}$ & $235.6$ \\
Mechanics-trained map & $156$ & $2.425\times 10^{-4}$ & $0.00504$ & $0.004911$ & $9.931\times 10^{-5}$ & $235.6$ \\
\bottomrule
\end{tabular}
\end{adjustbox}
\end{table}

\FloatBarrier

The six panels in Fig.~\ref{fig:vectorfields} connect the error measures in
Table~\ref{tab:equal} to their spatial origin.  The displacement, von Mises
stress, and strain-energy fields retain the expected mixed-mode pattern, while
the stress and energy concentrate at the crack tip.  Because the two
displacement-error panels use the same color scale, their contrast shows that
the learned coordinate suppresses the tip-dominated error rather than merely
rescaling its display.  The principal-stress directions remain consistent
with the physical mixed-mode field; the mapping improves radial resolution
without rotating the constitutive response or smoothing away the crack.

Fig.~\ref{fig:scalarfields} isolates the same mechanism in the scalar problem
used for training.  The gradient streamlines terminate at the singular point,
and the identity-map error follows this unresolved radial gradient.  The
mapped error is smaller on the common scale while the global solution and flux
directions are preserved.  This comparison is important for the subsequent
elasticity result: it shows that the transferred improvement originates in the
regularity of the coordinate learned from the scalar weak form, not in a
vector-field fit.

\paragraph{Sensitivity to initialization.}

\begin{figure}[htp]
\centering
\begin{subfigure}[t]{0.33\linewidth}\centering
\includegraphics[width=\linewidth]{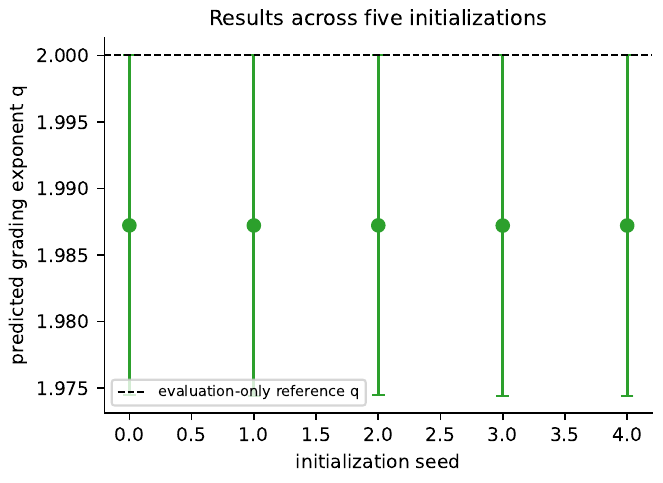}
\caption{Predicted radial grading exponents.}
\end{subfigure}\hfill
\begin{subfigure}[t]{0.33\linewidth}\centering
\includegraphics[width=\linewidth]{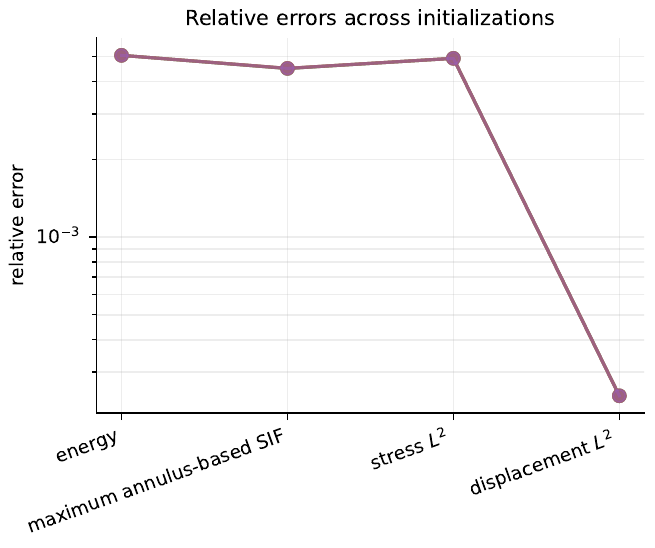}
\caption{Field-error variation across initializations.}
\end{subfigure}\hfill
\begin{subfigure}[t]{0.33\linewidth}\centering
\includegraphics[width=\linewidth]{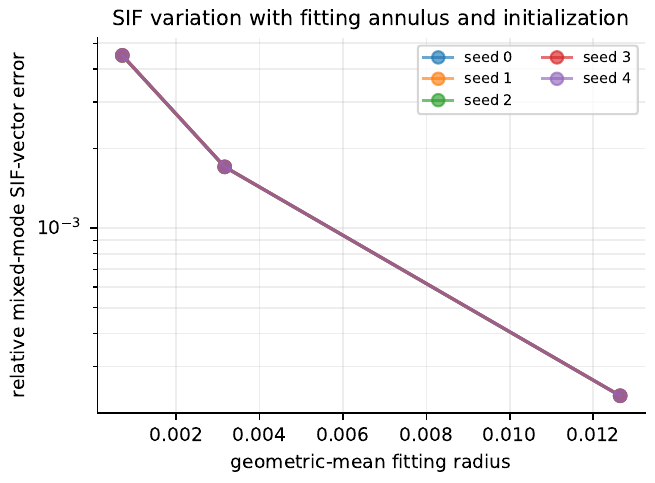}
\caption{Sensitivity to the fitting annulus.}
\end{subfigure}
\caption{Predicted radial grading exponents and cross-initialization field and
SIF errors.}
\label{fig:vectorseeds}
\end{figure}

Independent initializations determine whether the transferred exponent depends
materially on the optimization path.  Each model is trained from the same
mechanics energies and its predicted $q$ is passed unchanged to the mixed-mode
solve; exact exponents and vector errors remain evaluation quantities.  All
runs satisfy the prescribed tolerances, and the relative energy-norm error departs
from its mean \PrimaryVectorEnergyMean{} by a population standard deviation of
only \PrimaryVectorEnergyStd{}.  The nearly identical vector errors therefore
follow from a stable coordinate rather than a favorable initialization.
Table~\ref{tab:vectorseeds} reports the full statistics, with ``maximum''
denoting the sample maximum of each nonnegative error metric; the corresponding
variations are shown in Fig.~\ref{fig:vectorseeds}.

\begin{table}[htp]
\centering
\ResultTableFont
\caption{Statistics over independent conditional-exponent trainings and
vector-solver transfers.}
\label{tab:vectorseeds}
\begin{tabular}{lrrrr}
\toprule
Metric & mean & population std. & maximum & fixed limit \\
\midrule
Exponent absolute error & $1.278\times 10^{-2}$ & $5.53\times 10^{-6}$ & $1.279\times 10^{-2}$ & $2\times 10^{-2}$ \\
Relative $L^2$ error & $2.428\times 10^{-4}$ & $3.021\times 10^{-10}$ & $2.428\times 10^{-4}$ & $3.5\times 10^{-4}$ \\
Relative energy-norm error & $5.057\times 10^{-3}$ & $1.484\times 10^{-8}$ & $5.057\times 10^{-3}$ & $6\times 10^{-3}$ \\
Relative stress error & $4.923\times 10^{-3}$ & $1.04\times 10^{-8}$ & $4.923\times 10^{-3}$ & $6\times 10^{-3}$ \\
Maximum annulus-based SIF error & $4.5\times 10^{-3}$ & $1.95\times 10^{-6}$ & $4.503\times 10^{-3}$ & $6\times 10^{-3}$ \\
\bottomrule
\end{tabular}
\end{table}

\FloatBarrier

Fig.~\ref{fig:vectorseeds} shows why the small standard deviations in
Table~\ref{tab:vectorseeds} are mechanically relevant.  The predicted radial grading exponents
are nearly indistinguishable across initializations, the transferred field
errors collapse onto the same values, and the annulus-dependent SIF curves
retain the same shape for every seed.  The last observation separates
extraction-radius sensitivity from optimization sensitivity: the recovered
SIF estimate changes with the physical annulus, but not with the starting point of
the neural parameters.  The straight-crack conclusion therefore does not rely
on a favorable local minimum.

\paragraph{Interaction-integral evaluation.}

\begin{figure}[t]
\centering
\begin{subfigure}[t]{0.45\linewidth}\centering
\includegraphics[width=\linewidth]{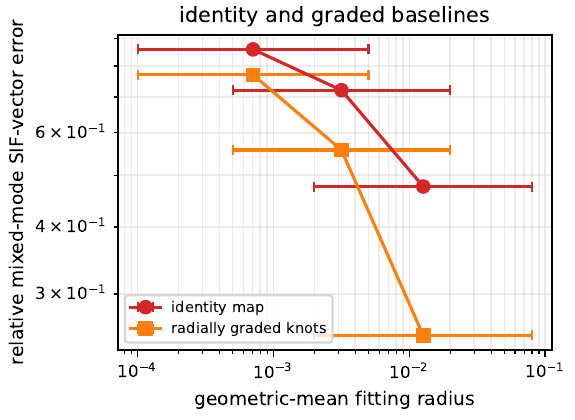}
\caption{Recovered stress intensity factors.}
\end{subfigure}\hfill
\begin{subfigure}[t]{0.45\linewidth}\centering
\includegraphics[width=\linewidth]{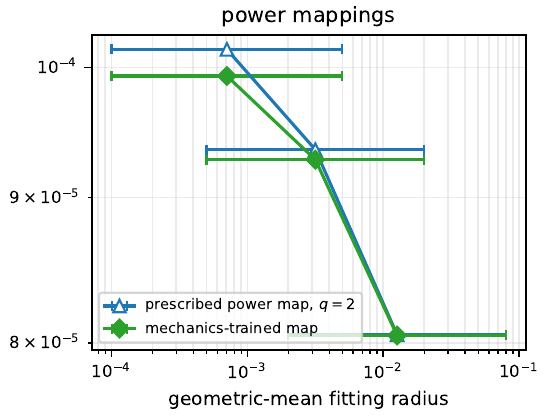}
\caption{Relative SIF error by fitting annulus.}
\end{subfigure}
\caption{Annulus-based least-squares estimates over physical fitting annuli.}
\label{fig:annuli}
\end{figure}

\begin{figure}[t]
\centering
\begin{subfigure}[t]{0.32\linewidth}\centering
\includegraphics[width=\linewidth]{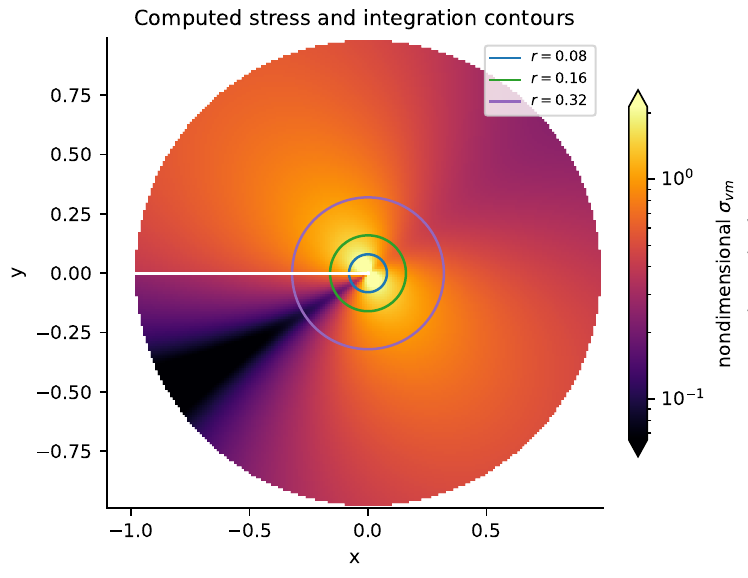}
\caption{Stress field and contours.}
\end{subfigure}\hfill
\begin{subfigure}[t]{0.32\linewidth}\centering
\includegraphics[width=\linewidth]{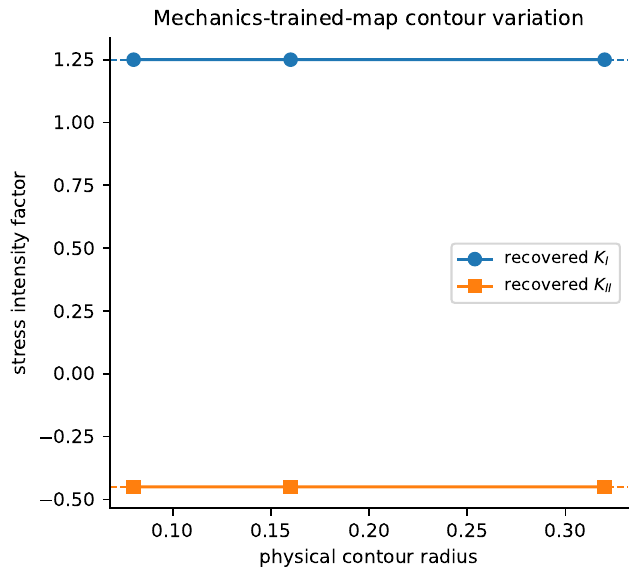}
\caption{Recovered stress intensity factors.}
\end{subfigure}\hfill
\begin{subfigure}[t]{0.32\linewidth}\centering
\includegraphics[width=\linewidth]{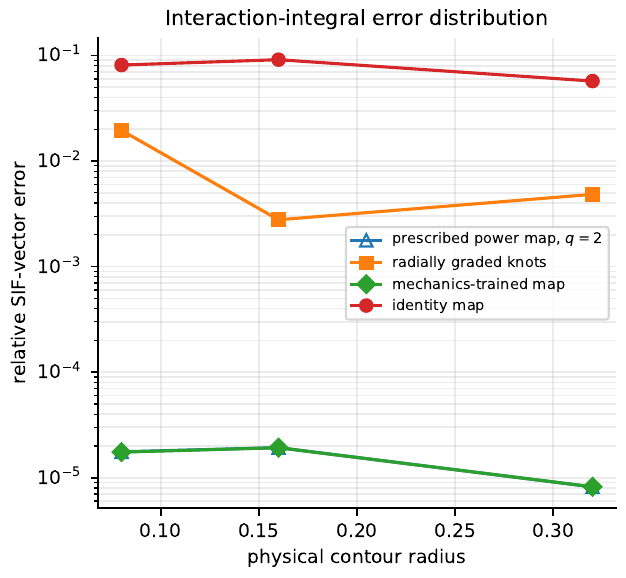}
\caption{Contour-wise relative error.}
\end{subfigure}
\caption{Stress field with interaction contours, recovered stress intensity
factors, and contour-wise relative errors.}
\label{fig:interaction}
\end{figure}

The interaction integral recovers the analytic mixed-mode field to
\InteractionExactMaxError{} and is insensitive to the bounded nonsingular $T$-stress
term.  For the learned map, the relative contour-to-contour variations are
\InteractionLearnedKISpread{} in $K_I$ and \InteractionLearnedKIISpread{} in
$K_{II}$, and the maximum error is \InteractionLearnedWorstError{}.  Its
agreement with the annulus-based displacement fit shows that the recovered
SIFs are not artifacts of the regression basis.  The same prescribed
contour radii are used for every coordinate in Table~\ref{tab:interaction};
Figs.~\ref{fig:annuli} and \ref{fig:interaction} show the annulus-based
displacement fit and interaction-integral results.

\begin{table}[t]
\centering
\ResultTableFont
\caption{Mixed-mode interaction-integral recovery on circular contours.}
\label{tab:interaction}
\begin{adjustbox}{max width=\linewidth}
\begin{tabular}{lrrrr}
\toprule
Coordinate & contour $r$ & recovered $K_I$ & recovered $K_{II}$ & rel. $e_K^{\mathcal I}$ \\
\midrule
Identity map & $0.08$ & $1.35$ & $-0.4893$ & $8.115\times 10^{-2}$ \\
 & $0.16$ & $1.364$ & $-0.4925$ & $9.125\times 10^{-2}$ \\
 & $0.32$ & $1.178$ & $-0.4264$ & $5.723\times 10^{-2}$ \\
\addlinespace
Radially graded knots & $0.08$ & $1.227$ & $-0.4371$ & $1.952\times 10^{-2}$ \\
 & $0.16$ & $1.247$ & $-0.4521$ & $2.777\times 10^{-3}$ \\
 & $0.32$ & $1.256$ & $-0.4533$ & $4.834\times 10^{-3}$ \\
\addlinespace
Prescribed power map $q=2$ & $0.08$ & $1.25$ & $-0.45$ & $1.759\times 10^{-5}$ \\
 & $0.16$ & $1.25$ & $-0.45$ & $1.926\times 10^{-5}$ \\
 & $0.32$ & $1.25$ & $-0.45$ & $8.181\times 10^{-6}$ \\
\addlinespace
Mechanics-trained map & $0.08$ & $1.25$ & $-0.45$ & $1.746\times 10^{-5}$ \\
 & $0.16$ & $1.25$ & $-0.45$ & $1.923\times 10^{-5}$ \\
 & $0.32$ & $1.25$ & $-0.45$ & $8.197\times 10^{-6}$ \\
\bottomrule
\end{tabular}
\end{adjustbox}
\end{table}

\FloatBarrier

The two panels of Fig.~\ref{fig:annuli} separate the large-error baselines from
the regularized coordinates.  The identity-map and radially graded-knot solutions retain a
marked dependence on annulus radius, whereas the prescribed-power and mechanics-trained
maps lie together at a much lower error level.  Their overlap is a substantive
result: for the known square-root crack field, learning correctly reduces to
selecting the classical power instead of inventing an unnecessary density
correction.

Fig.~\ref{fig:interaction} tests the same conclusion with a different
fracture functional.  The circular contours surround the high-stress region,
yet the recovered $K_I$ and $K_{II}$ remain close to their prescribed values as
the contour radius changes; the small contour-wise errors quantify this
numerical contour invariance.  Agreement between the annulus-based displacement fit and the
interaction-integral results is important because neither functional participates in training.
It therefore rules out the possibility that the improvement is specific to
one post-processing basis or one extraction support.

\paragraph{Refinement and conditioning.}

\begin{figure}[t]
\centering
\begin{subfigure}[t]{0.45\linewidth}\centering
\includegraphics[width=\linewidth]{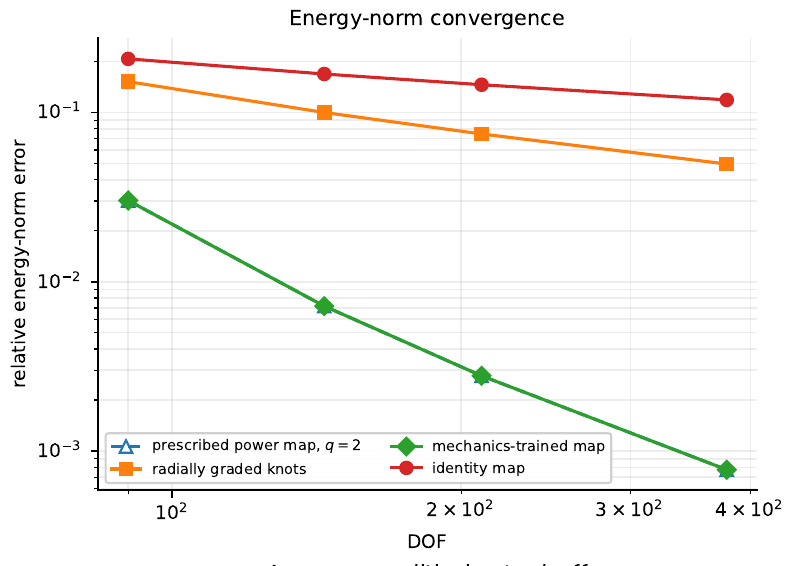}
\caption{Energy-error convergence.}
\end{subfigure}\hfill
\begin{subfigure}[t]{0.45\linewidth}\centering
\includegraphics[width=\linewidth]{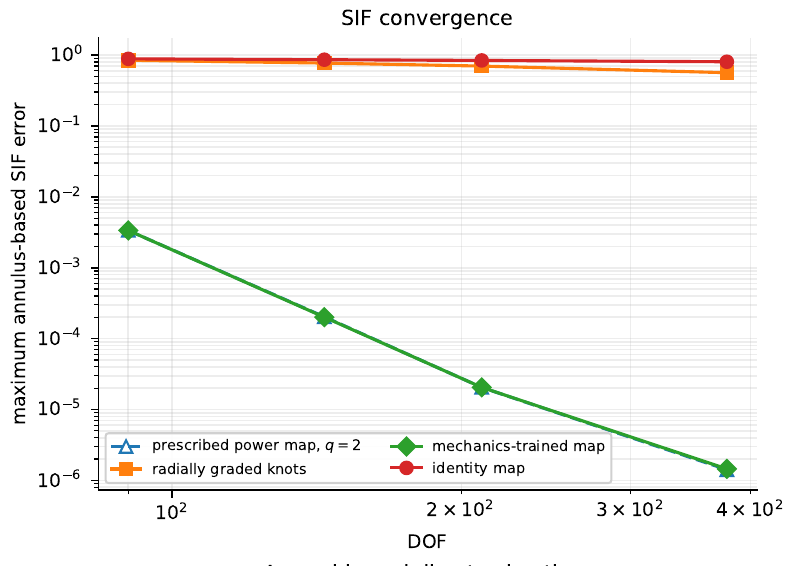}
\caption{Annulus-based SIF-error convergence.}
\end{subfigure}

\medskip
\begin{subfigure}[t]{0.45\linewidth}\centering
\includegraphics[width=\linewidth]{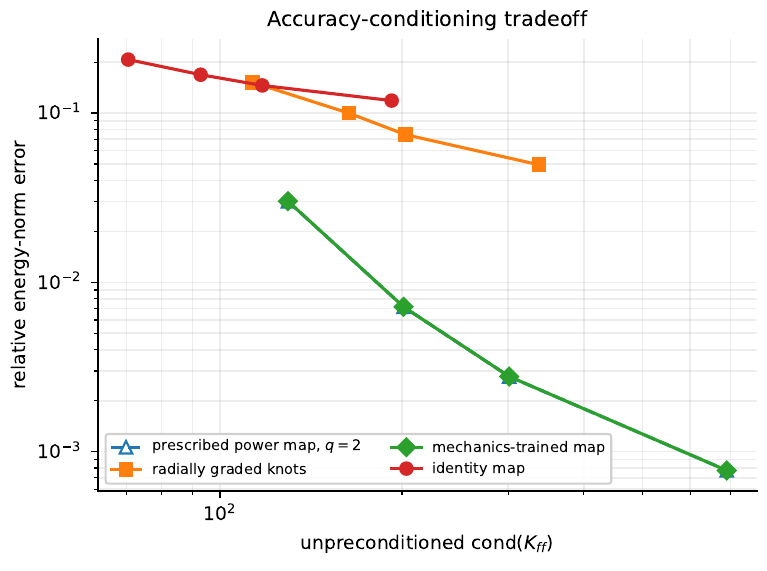}
\caption{Accuracy--conditioning trade-off.}
\end{subfigure}\hfill
\begin{subfigure}[t]{0.45\linewidth}\centering
\includegraphics[width=\linewidth]{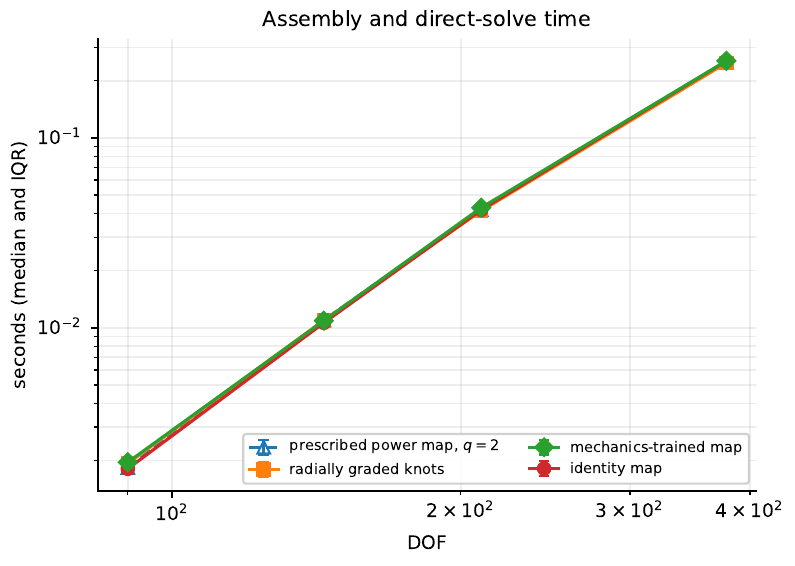}
\caption{Assembly and solution time.}
\end{subfigure}
\caption{Energy-norm and SIF convergence, conditioning, and
assembly/solution time.}
\label{fig:convergence}
\end{figure}

\begin{figure}[t]
\centering
\begin{subfigure}[t]{0.32\linewidth}\centering
\includegraphics[width=\linewidth]{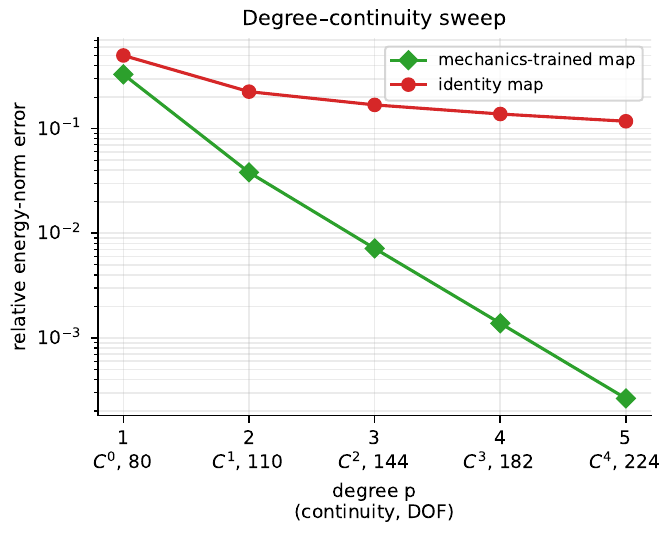}
\caption{Displacement error.}
\end{subfigure}\hfill
\begin{subfigure}[t]{0.32\linewidth}\centering
\includegraphics[width=\linewidth]{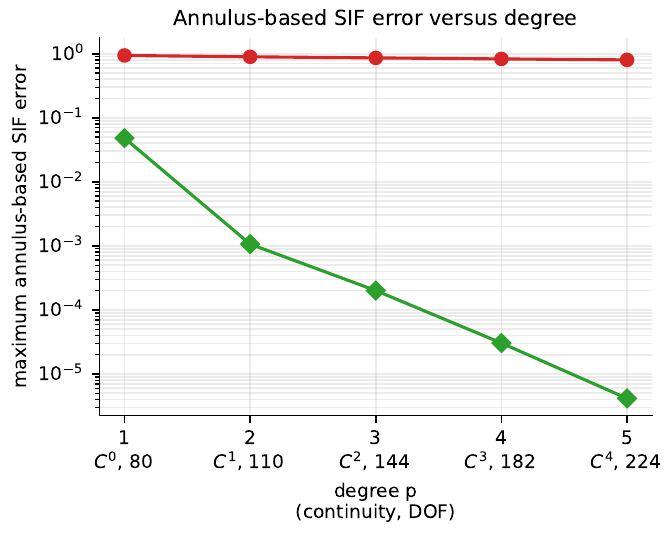}
\caption{Energy error.}
\end{subfigure}\hfill
\begin{subfigure}[t]{0.32\linewidth}\centering
\includegraphics[width=\linewidth]{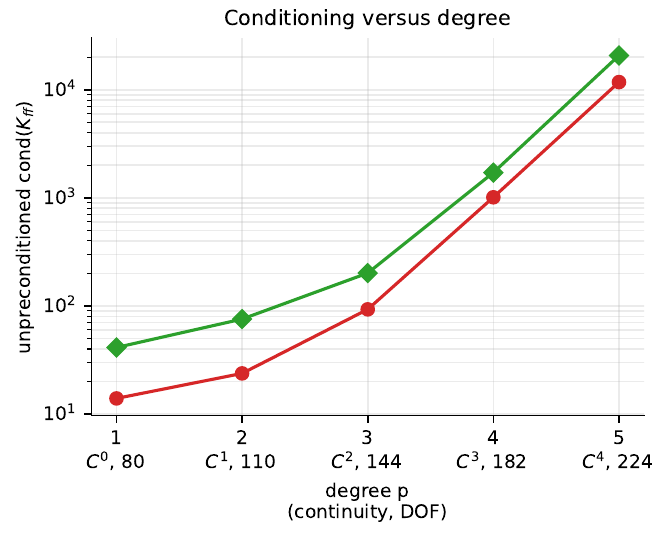}
\caption{Maximum annulus-based SIF error.}
\end{subfigure}
\caption{Influence of spline degree and continuity on crack-solution errors.}
\label{fig:degree}
\end{figure}

Over the selected fine-grid interval, the mechanics-trained and prescribed-power maps attain a log--log slope of energy-norm error versus DOF of $-2.295$, compared with $-0.3636$ and $-0.7238$ for the
identity and radially graded mappings.  The improvement is stronger for the annulus-based
SIF error, as Table~\ref{tab:slopes} and Fig.~\ref{fig:convergence} show.
Fig.~\ref{fig:degree} confirms that the same ordering persists when the spline
degree and continuity are changed.  The timing curves in
Fig.~\ref{fig:convergence} nearly coincide at the plotted scale; they comprise
assembly and direct solution rather than the end-to-end cost considered below.

\begin{table}[htp]
\centering
\ResultTableFont
\caption{Fitted fine-grid slopes from the refinement sweep.}
\label{tab:slopes}
\begin{tabular}{lrrr}
\toprule
Coordinate & fitted DOF interval & energy--DOF slope & SIF--DOF slope \\
\midrule
Identity map & $144-378$ & $-0.3636$ & $-0.06937$ \\
Radially graded knots & $144-378$ & $-0.7238$ & $-0.3277$ \\
Prescribed power map $q=2$ & $144-378$ & $-2.295$ & $-5.094$ \\
Mechanics-trained map & $144-378$ & $-2.295$ & $-5.054$ \\
\bottomrule
\end{tabular}
\end{table}

\FloatBarrier

In Fig.~\ref{fig:convergence}, the steeper energy and SIF curves of the fixed
and learned maps show that their advantage grows under refinement rather than
being confined to one matched-size comparison.  The conditioning panel records
the price of concentrating resolution near the collapsed edge, but the timing
panel shows that the four assembly-and-solve curves remain nearly coincident at
equal dimension.  Thus the observed accuracy gain is not produced by hidden
extra linear-solver work.  Fig.~\ref{fig:degree} reaches the same ordering
when degree and continuity are varied.  This robustness is important because
it identifies the coordinate regularity, rather than one particular spline
space, as the controlling factor; the overlap of the prescribed-power and mechanics-trained curves
also anticipates the numerically negligible density weights reported below.

The optimizer selects $q=\LearnedQ{}$, but the two trained density weights are
\LearnedWeightOne{} and \LearnedWeightTwo{}.  The mechanics-trained-to-prescribed-power energy-norm error
ratio is \FixedLearnedEnergyRatio{}, and the corresponding annulus-based SIF error
ratio is \FixedLearnedSIFRatio{}.  The straight-crack improvement is therefore
produced by the selected grading exponent: the two density weights are numerically negligible, and
the prescribed power map reproduces the gain.  This comparison separates
radial grading-exponent selection from a genuinely non-power-law neural correction.

\subsection{Conditional wedge problems}

\begin{figure}[htp]
\centering
\begin{subfigure}[t]{0.32\linewidth}\centering
\includegraphics[width=\linewidth]{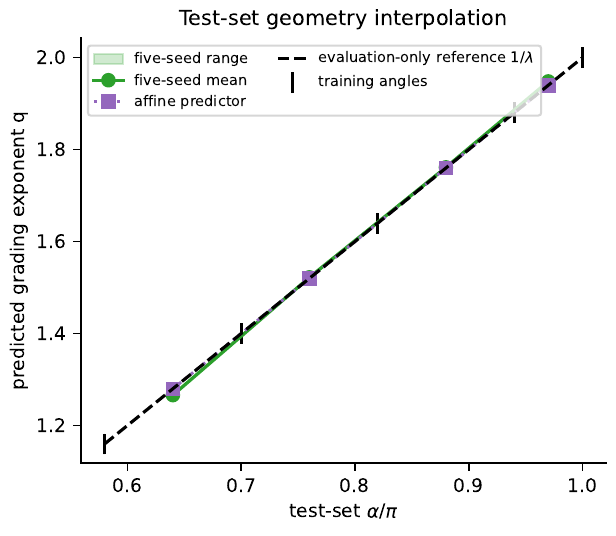}
\caption{Conditional prediction.}
\end{subfigure}\hfill
\begin{subfigure}[t]{0.32\linewidth}\centering
\includegraphics[width=\linewidth]{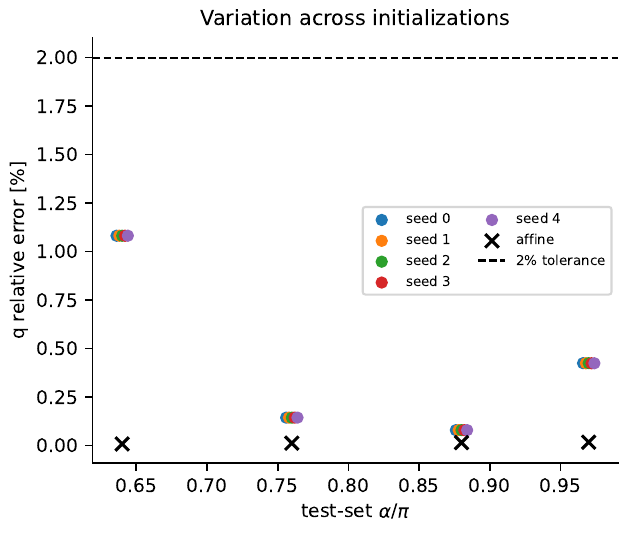}
\caption{Test-set grading-exponent error.}
\end{subfigure}\hfill
\begin{subfigure}[t]{0.32\linewidth}\centering
\includegraphics[width=\linewidth]{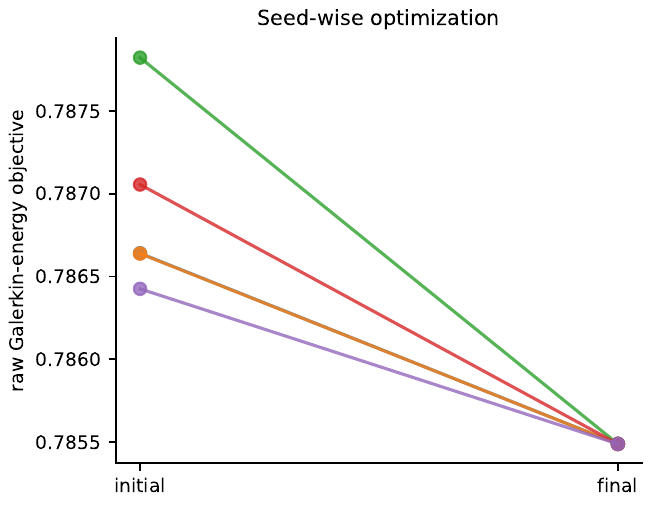}
\caption{Initialization distribution.}
\end{subfigure}
\caption{Wedge test-set predictions and their variation across
initializations.}
\label{fig:wedge}
\end{figure}

\begin{figure}[htp]
\centering
\begin{subfigure}[t]{0.45\linewidth}\centering
\includegraphics[width=\linewidth]{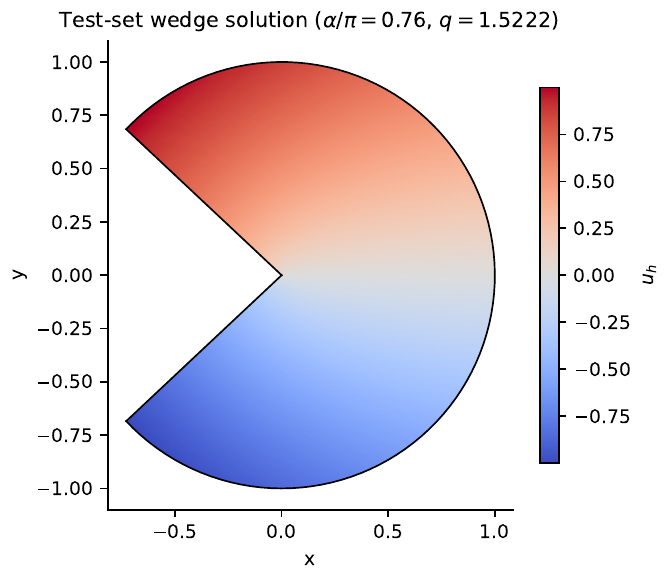}
\caption{Solved field and gradient.}
\end{subfigure}\hfill
\begin{subfigure}[t]{0.45\linewidth}\centering
\includegraphics[width=\linewidth]{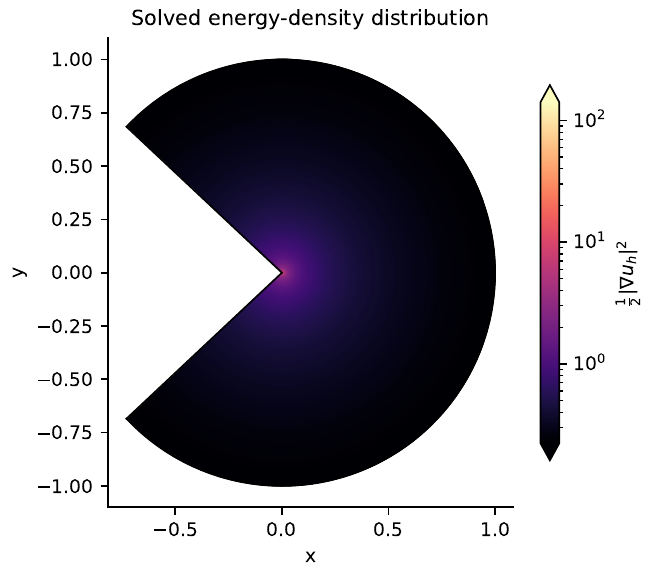}
\caption{Local energy density.}
\end{subfigure}
\caption{Computed wedge field, gradient, and local energy distributions at a
test geometry.}
\label{fig:wedgefield}
\end{figure}

The affine predictor attains the smaller test-set grading-exponent error because the
manufactured angle-to-grading-exponent relation is itself affine.  Its maximum
error is \WedgeAffineWorstQError{}, compared with \WedgeANNWorstQError{} for the
neural rule, although both are trained from the same Galerkin energy without
exponent labels.  The neural parameterization therefore adds no useful
approximation capacity to this family.  Table~\ref{tab:wedge} gives the
test-set errors, while Figs.~\ref{fig:wedge} and \ref{fig:wedgefield} show their
variation and a representative wedge field.

\begin{table}[htp]
\centering
\ResultTableFont
\caption{Interpolation on the wedge test set.}
\label{tab:wedge}
\begin{tabular}{lrrr}
\toprule
Conditional exponent model & mean rel. $q$ error & max. rel. $q$ error & max. rel. energy \\
\midrule
Bounded neural predictor & $4.32\times 10^{-3}$ & $1.082\times 10^{-2}$ & $5.683\times 10^{-4}$ \\
Affine predictor & $1.285\times 10^{-4}$ & $1.667\times 10^{-4}$ & $1.194\times 10^{-4}$ \\
\bottomrule
\end{tabular}
\end{table}

\FloatBarrier

Fig.~\ref{fig:wedge} makes the affine family a limiting case for assessing neural
flexibility.  The affine predictions follow the manufactured
angle-to-grading-exponent relation more closely at the test angles, while the
neural predictions show a larger systematic interpolation error despite their
small variation across initializations.  The comparison demonstrates that a
network is not advantageous merely because it can represent a more complicated
rule.  In Fig.~\ref{fig:wedgefield}, the gradient and local energy concentrate
at the re-entrant corner selected by that exponent.  The field plot therefore
links the scalar prediction error to the solved Galerkin mechanics: an error in
$q$ changes the resolution of the corner field, not just an auxiliary label.

\paragraph{Nonlinear extrapolation.}

Unlike the preceding wedge, the Robin family combines a nonlinear
parameter-to-grading-exponent relation with two-sided extrapolation.  Its training
parameters lie in $[\NonlinearTrainingMinimum{},\NonlinearTrainingMaximum{}]$,
whereas the two test parameters lie beyond opposite ends of this interval.  The
post-training affine-fit residual, \NonlinearAffineResidual{}, confirms that an
affine conditional rule cannot represent the evaluated exponent family.

\begin{figure}[t]
\centering
\begin{subfigure}[t]{0.32\linewidth}\centering
\includegraphics[width=\linewidth]{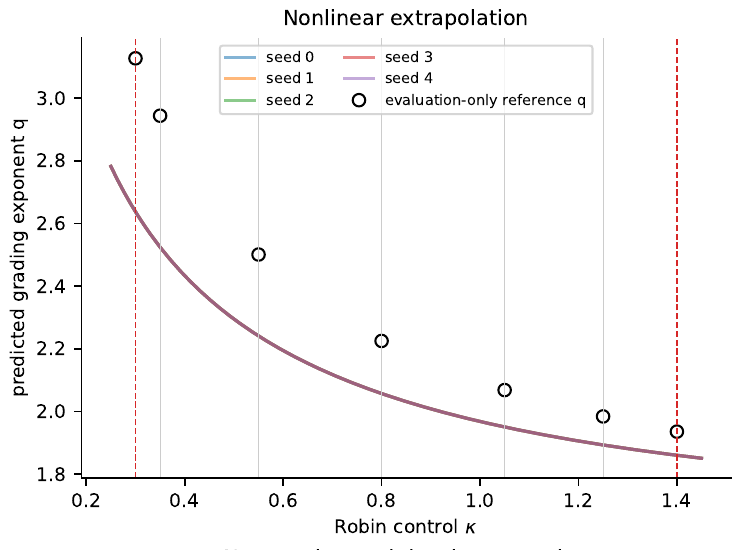}
\caption{Nonlinear grading-exponent prediction.}
\end{subfigure}\hfill
\begin{subfigure}[t]{0.32\linewidth}\centering
\includegraphics[width=\linewidth]{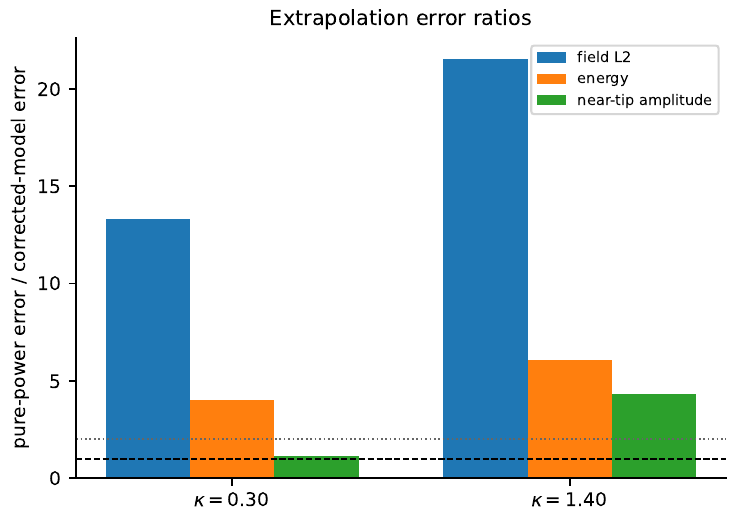}
\caption{Extrapolation error ratios.}
\end{subfigure}\hfill
\begin{subfigure}[t]{0.32\linewidth}\centering
\includegraphics[width=\linewidth]{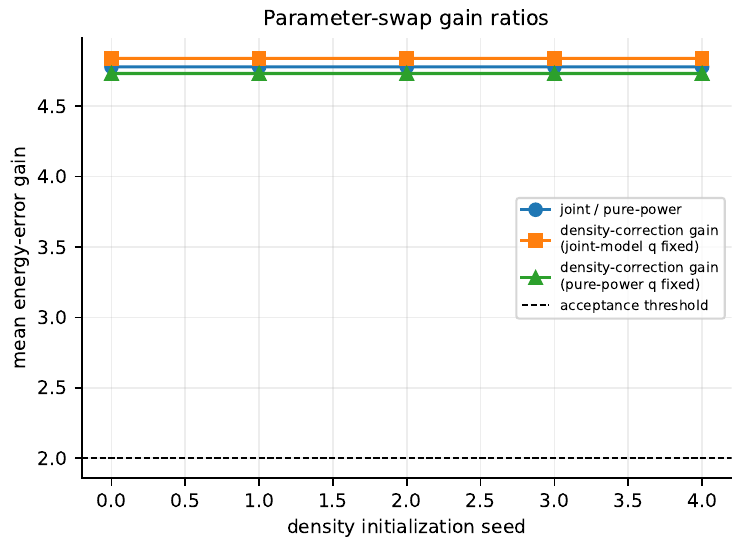}
\caption{Parameter-swap gain ratios.}
\end{subfigure}

\medskip
\begin{subfigure}[t]{0.32\linewidth}\centering
\includegraphics[width=\linewidth]{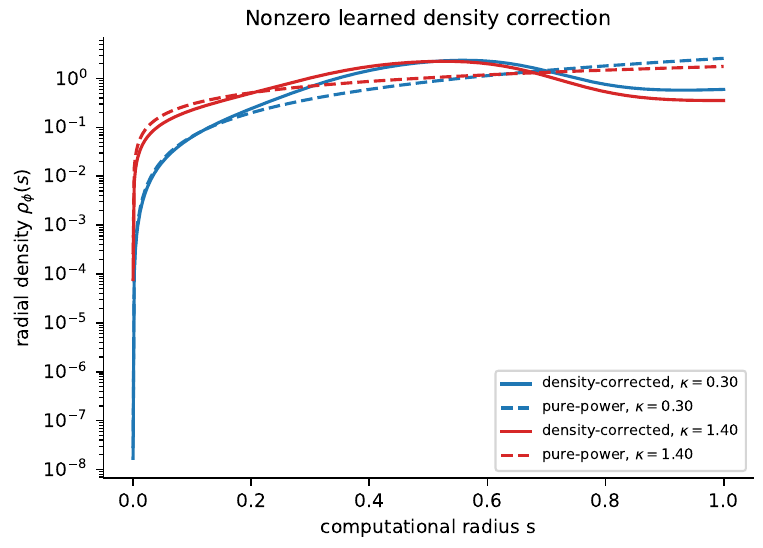}
\caption{Learned radial-density correction.}
\end{subfigure}\hfill
\begin{subfigure}[t]{0.32\linewidth}\centering
\includegraphics[width=\linewidth]{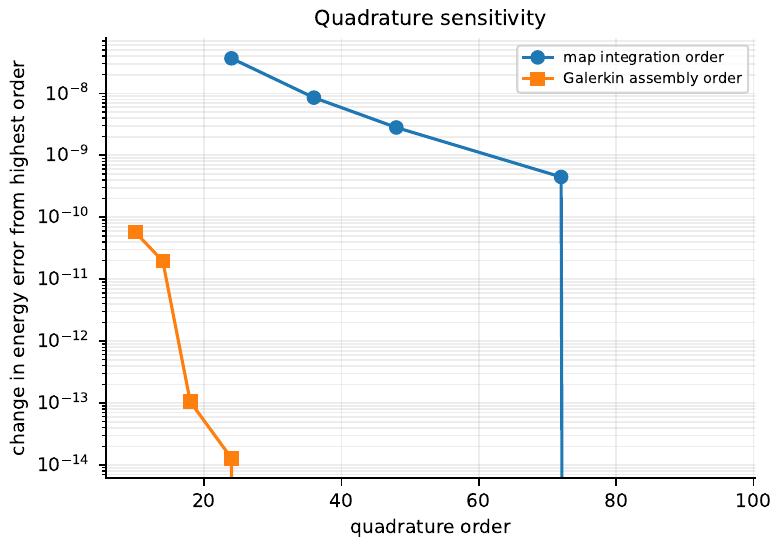}
\caption{Quadrature sensitivity.}
\end{subfigure}\hfill
\begin{subfigure}[t]{0.32\linewidth}\centering
\includegraphics[width=\linewidth]{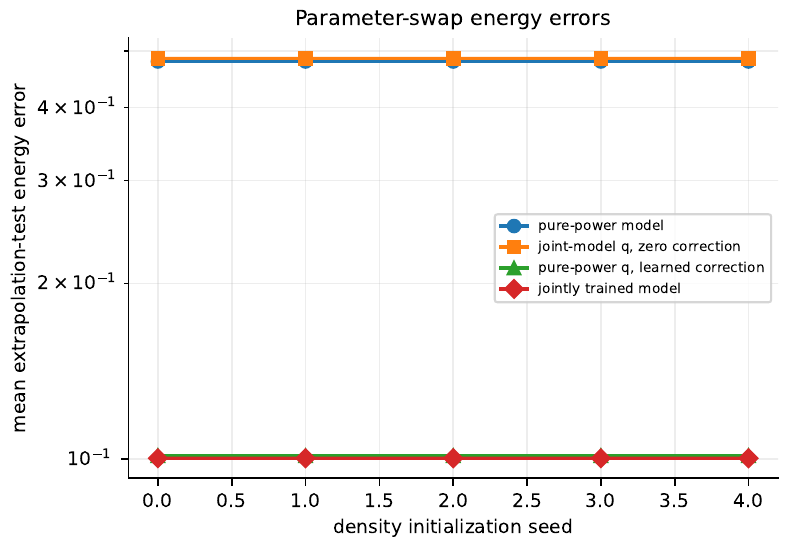}
\caption{Parameter-swap energy errors.}
\end{subfigure}
\caption{Nonlinear extrapolation, density correction, ablation, and quadrature
results.}
\label{fig:nonlineartransfer}
\end{figure}

\begin{figure}[htp]
\centering
\begin{subfigure}[t]{0.45\linewidth}\centering
\includegraphics[width=\linewidth]{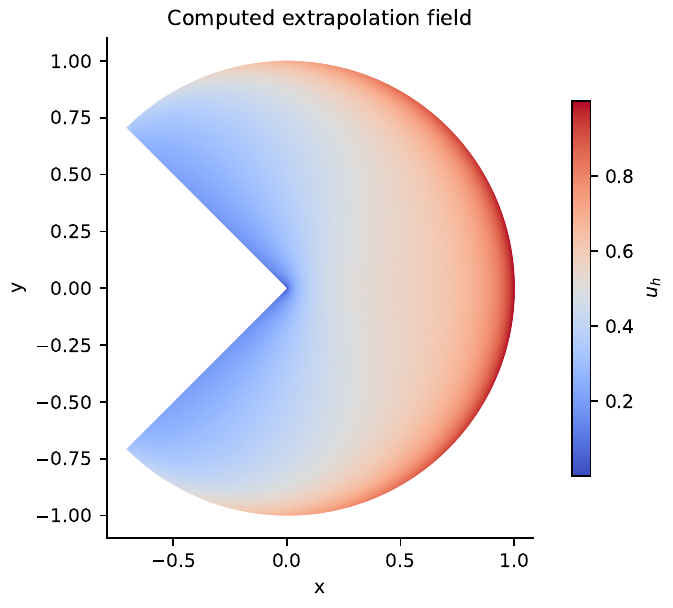}
\caption{Solved field distribution.}
\end{subfigure}\hfill
\begin{subfigure}[t]{0.45\linewidth}\centering
\includegraphics[width=\linewidth]{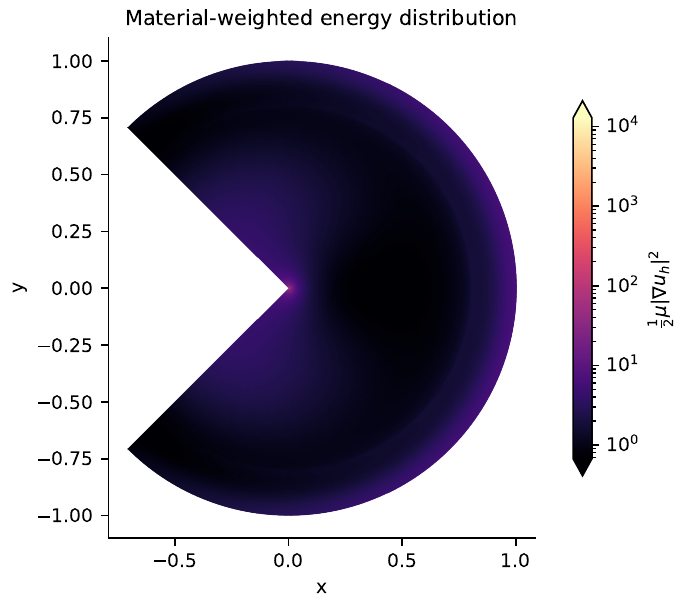}
\caption{Material-weighted energy.}
\end{subfigure}

\medskip
\begin{subfigure}[t]{0.45\linewidth}\centering
\includegraphics[width=\linewidth]{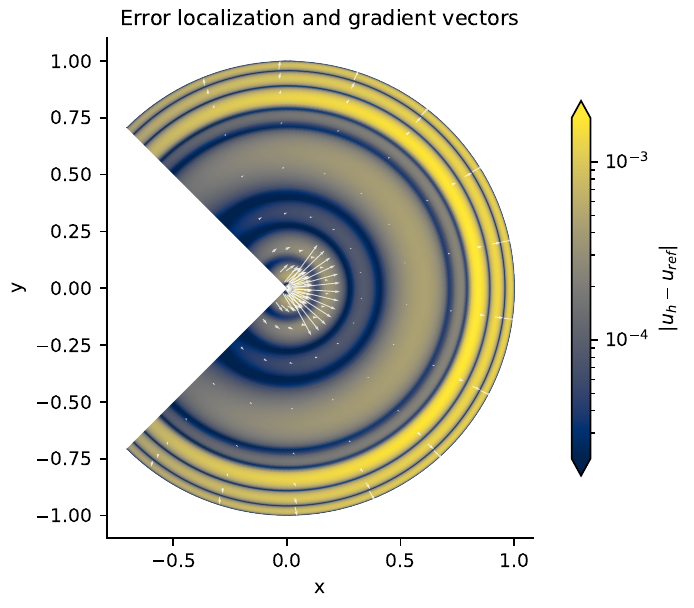}
\caption{Localized error and gradient.}
\end{subfigure}\hfill
\begin{subfigure}[t]{0.45\linewidth}\centering
\includegraphics[width=\linewidth]{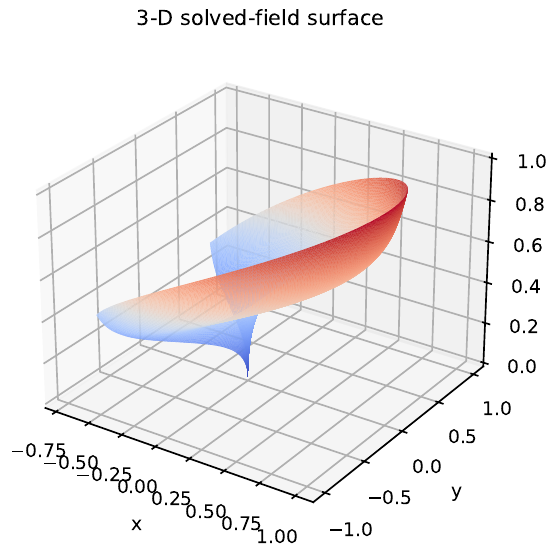}
\caption{Three-dimensional field surface.}
\end{subfigure}
\caption{Field and error distributions at the upper extrapolation test parameter.}
\label{fig:nonlinearfield}
\end{figure}

Every accepted initialization retains a nonzero density correction.  Relative
to the shared pure-power model, it improves the field, energy, and near-tip
amplitude errors by at least \NonlinearFieldGainMinimum{},
\NonlinearEnergyGainMinimum{}, and \NonlinearAmplitudeGainMinimum{},
respectively.  The improvement therefore persists across both extrapolation test parameters and all accepted initializations.  In Table~\ref{tab:nonlinear}, the
gain is $e_E^{\rm power}/e_E^{\rm corrected}$, so values greater than one favor the
density-corrected map.  Figs.~\ref{fig:nonlineartransfer} and
\ref{fig:nonlinearfield} show the learned rule and a representative
extrapolation field.

\begin{table}[t]
\centering
\ResultTableFont
\caption{Nonlinear-family extrapolation comparison.}
\label{tab:nonlinear}
\begin{adjustbox}{max width=\linewidth}
\begin{tabular}{lrrrrr}
\toprule
Extrapolation parameter & pure-power $e_E$ & corrected mean $e_E$ & corrected std. & corrected max. & energy gain \\
\midrule
$0.3$ & $0.5025$ & $0.1253$ & $8.179\times 10^{-7}$ & $0.1253$ & $4.012$ \\
$1.4$ & $0.4569$ & $0.07549$ & $4.07\times 10^{-7}$ & $0.07549$ & $6.052$ \\
\bottomrule
\end{tabular}
\end{adjustbox}
\end{table}

The parameter-swap ablation uses $P=(q_P,\bm0)$ for the pure-power model,
$A=(q_J,\bm0)$, $B=(q_P,\bm w_J)$, and
$J=(q_J,\bm w_J)$ for the jointly trained density-corrected model; none of
these combinations is refitted.  Removing $\bm w$ from $J$
increases the energy error by at least
\NonlinearDensityIncrementalEnergyGainMinimum{}, while adding it to the
pure-power grading exponent gives a gain of
\NonlinearDensityOnPowerQEnergyGainMinimum{}.  The non-power-law correction thus
contributes independently of the exponent predictor.  Table~\ref{tab:nonlinearablation}
reports the mean relative energy errors; ratios greater than one favor the
denominator, and none of the parameter combinations is retrained.

\begin{table}[t]
\centering
\ResultTableFont
\caption{Parameter-swap ablation at the two extrapolation test parameters.}
\label{tab:nonlinearablation}
\begin{adjustbox}{max width=\linewidth}
\begin{tabular}{rrrrrrrr}
\toprule
Seed & $e_E^P$ & $e_E^A$ & $e_E^B$ & $e_E^J$ & $P/J$ & $A/J$ & $P/B$ \\
\midrule
$0$ & $0.4797$ & $0.4858$ & $0.1014$ & $0.1004$ & $4.779$ & $4.839$ & $4.732$ \\
$1$ & $0.4797$ & $0.4858$ & $0.1014$ & $0.1004$ & $4.779$ & $4.839$ & $4.732$ \\
$2$ & $0.4797$ & $0.4858$ & $0.1014$ & $0.1004$ & $4.779$ & $4.839$ & $4.732$ \\
$3$ & $0.4797$ & $0.4858$ & $0.1014$ & $0.1004$ & $4.779$ & $4.839$ & $4.732$ \\
$4$ & $0.4797$ & $0.4858$ & $0.1014$ & $0.1004$ & $4.779$ & $4.839$ & $4.732$ \\
\bottomrule
\end{tabular}
\end{adjustbox}
\end{table}

Changing the radial-map normalization and Galerkin quadrature orders changes the
energy error by at most \MapQuadratureEnergyDelta{} and
\AssemblyQuadratureEnergyDelta{}, respectively.  These variations are smaller
than the density-correction gain, so the correction is not a quadrature-order effect.
The deltas in Table~\ref{tab:quadrature} are measured from the corresponding
highest-order result.

\begin{table}[t]
\centering
\ResultTableFont
\caption{Quadrature sensitivity of the density-corrected nonlinear map.}
\label{tab:quadrature}
\begin{adjustbox}{max width=\linewidth}
\begin{tabular}{lrrrrrr}
\toprule
Sweep & map order & assembly order & $\Delta e_{L^2}$ & $\Delta e_E$ & $\Delta e_A$ & min. density \\
\midrule
Map normalization & $24$ & $32$ & $1.225\times 10^{-7}$ & $3.681\times 10^{-8}$ & $5.871\times 10^{-7}$ & $1.475\times 10^{-6}$ \\
 & $36$ & $32$ & $2.71\times 10^{-8}$ & $8.542\times 10^{-9}$ & $1.302\times 10^{-7}$ & $1.475\times 10^{-6}$ \\
 & $48$ & $32$ & $8.924\times 10^{-9}$ & $2.812\times 10^{-9}$ & $4.249\times 10^{-8}$ & $1.475\times 10^{-6}$ \\
 & $72$ & $32$ & $1.419\times 10^{-9}$ & $4.473\times 10^{-10}$ & $7.243\times 10^{-9}$ & $1.475\times 10^{-6}$ \\
 & $96$ & $32$ & $0$ & $0$ & $0$ & $1.475\times 10^{-6}$ \\
\addlinespace
Galerkin assembly & $96$ & $10$ & $1.133\times 10^{-8}$ & $5.738\times 10^{-11}$ & $7.472\times 10^{-8}$ & $1.475\times 10^{-6}$ \\
 & $96$ & $14$ & $3.848\times 10^{-10}$ & $1.955\times 10^{-11}$ & $9.562\times 10^{-10}$ & $1.475\times 10^{-6}$ \\
 & $96$ & $18$ & $1.987\times 10^{-12}$ & $1.051\times 10^{-13}$ & $5.376\times 10^{-11}$ & $1.475\times 10^{-6}$ \\
 & $96$ & $24$ & $2.481\times 10^{-13}$ & $1.285\times 10^{-14}$ & $5.592\times 10^{-12}$ & $1.475\times 10^{-6}$ \\
 & $96$ & $32$ & $0$ & $0$ & $0$ & $1.475\times 10^{-6}$ \\
\bottomrule
\end{tabular}
\end{adjustbox}
\end{table}

The representative fields in Fig.~\ref{fig:nonlinearfield} correspond to the
upper extrapolation test parameter $\kappa=\NonlinearRepresentativeControl{}$ and to
the admissible initialization with the smallest density-corrected-model objective.

\FloatBarrier

Fig.~\ref{fig:nonlineartransfer} separates three explanations for the
extrapolation improvement.  The predicted grading-exponent curves capture the overall
nonlinear parameter dependence, but the independent field, energy, and amplitude
gains remain above unity only when the radial-density correction is nonzero.
The parameter swaps show that this gain cannot be assigned to $q$ alone,
and the quadrature panels show that it is larger than the numerical variation
caused by either integration order.  The agreement across seeds further shows
that the correction is a property of the learned coordinate family rather than
one initialization.

The spatial evidence in Fig.~\ref{fig:nonlinearfield} explains why the extra
density freedom matters.  The global solution remains smooth away from the
singular point, whereas the material-weighted energy, gradient, and residual
error concentrate in the near-tip region whose radial sampling is changed by
$\bm w$.  The three-dimensional view exposes the same radial structure without
introducing a second model.  These panels are therefore retained as a
mechanistic check: the non-power-law correction acts where the singular weak form
is most sensitive, rather than lowering an aggregate loss through a remote
change in the field.

\subsection{Map admissibility and computational cost}

\begin{figure}[htp]
\centering
\begin{subfigure}[t]{0.45\linewidth}\centering
\includegraphics[width=\linewidth]{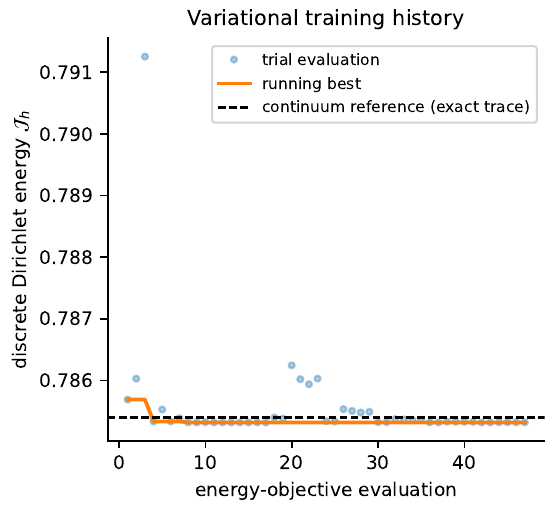}
\caption{Equilibrium energy-objective history.}
\end{subfigure}\hfill
\begin{subfigure}[t]{0.45\linewidth}\centering
\includegraphics[width=\linewidth]{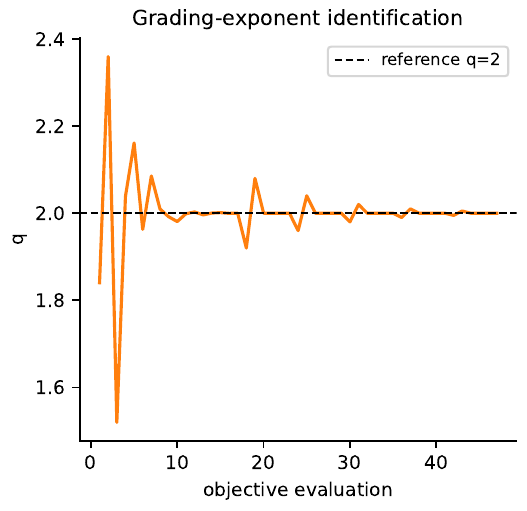}
\caption{Learned radial coordinate.}
\end{subfigure}

\medskip
\begin{subfigure}[t]{0.45\linewidth}\centering
\includegraphics[width=\linewidth]{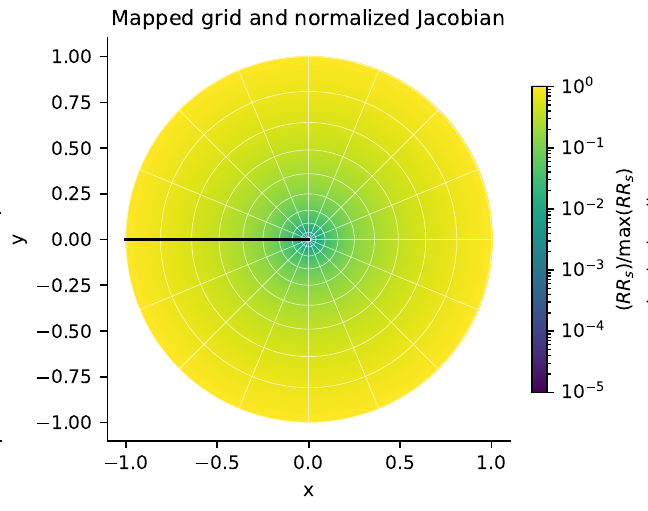}
\caption{Positive coordinate derivative.}
\end{subfigure}\hfill
\begin{subfigure}[t]{0.45\linewidth}\centering
\includegraphics[width=\linewidth]{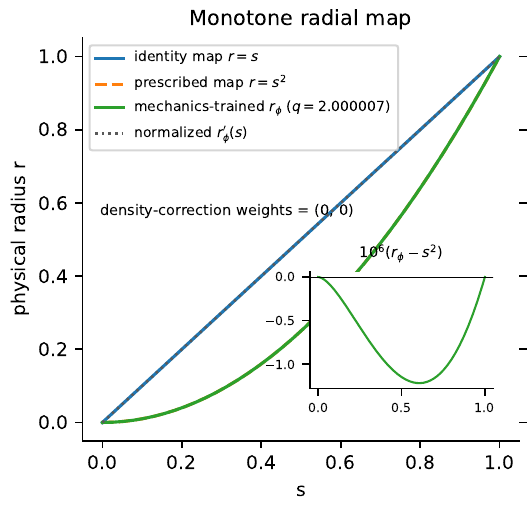}
\caption{Mapped coordinate grid.}
\end{subfigure}
\caption{Training history and the resulting radial coordinate and grid.}
\label{fig:maps}
\end{figure}
\FloatBarrier

Fig.~\ref{fig:maps} links the optimization history to the coordinate that is
actually used by the Galerkin solver.  A decreasing equilibrium energy objective alone
would not establish a valid map; the radial graph, the strictly positive
derivative, and the non-crossing mapped grid provide the complementary
geometric evidence.  The grid also shows where resolution has been moved:
radial lines remain ordered while cells contract toward the singular point.
This figure is essential to the accuracy results because an error reduction
obtained from a folded or orientation-reversing chart would have no numerical
meaning.

\paragraph{Jacobian tests.}

\begin{figure}[t]
\centering
\begin{subfigure}[t]{0.45\linewidth}\centering
\includegraphics[width=\linewidth]{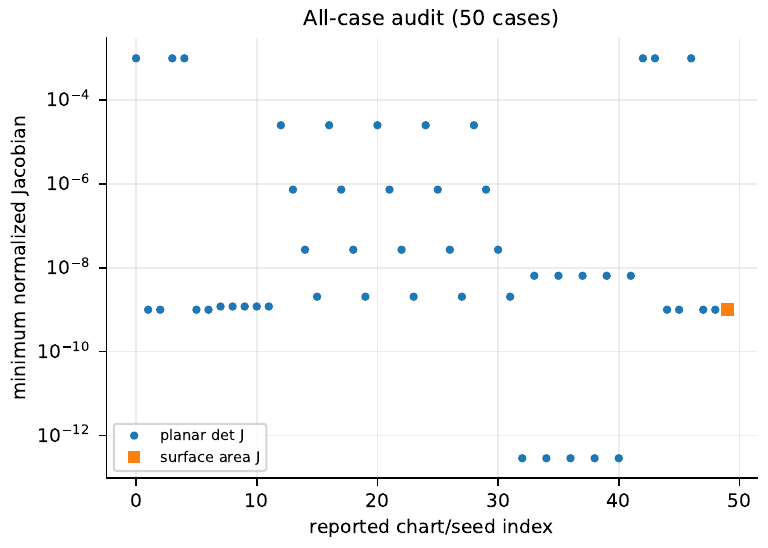}
\caption{Case-wise Jacobian minima.}
\end{subfigure}\hfill
\begin{subfigure}[t]{0.45\linewidth}\centering
\includegraphics[width=\linewidth]{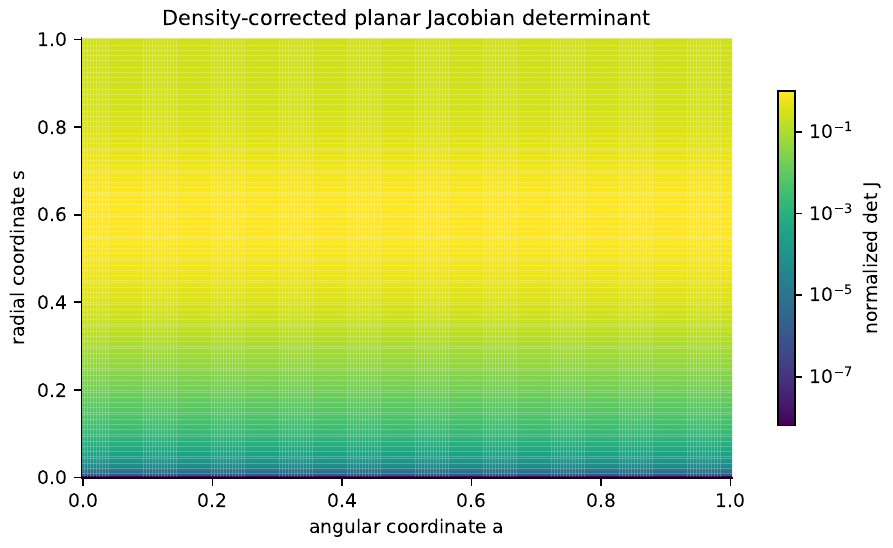}
\caption{Density-corrected nonlinear planar chart.}
\end{subfigure}

\caption{Case-wise normalized Jacobian minima and the density-corrected nonlinear planar Jacobian determinant.}
\label{fig:jacobianaudit}
\end{figure}

Away from the collapsed edge, an admissible chart must preserve orientation
and its implemented derivatives must agree with an independent approximation.
For the common sample grid $\mathcal G$, restricted to the collapsed-edge
cutoff $s\ge\JacobianCutoff{}$, let $J$ be the
signed planar Jacobian determinant relative to the reference orientation or the surface Jacobian.
We measure these two requirements by
\begin{align}
 \widehat J&=J/\max_{\mathcal G}J,\nonumber\\
 \delta_D&=\max\!\left\{
 \frac{\max_{\mathcal G}\norm{\bm\chi_{\phi,s}^{\mathrm{an}}-\bm\chi_{\phi,s}^{\mathrm{FD}}}}
      {\max_{\mathcal G}\norm{\bm\chi_{\phi,s}^{\mathrm{an}}}},
 \frac{\max_{\mathcal G}\norm{\bm\chi_{\phi,a}^{\mathrm{an}}-\bm\chi_{\phi,a}^{\mathrm{FD}}}}
      {\max_{\mathcal G}\norm{\bm\chi_{\phi,a}^{\mathrm{an}}}}
 \right\}.
 \label{eq:jacobianaudit}
\end{align}
In Eq.~\eqref{eq:jacobianaudit}, ``an'' and ``FD'' denote the analytic and
finite-difference derivatives.
All \JacobianCaseCount{} charts preserve the sampled orientation, with minimum
normalized Jacobian \JacobianMinimumNormalized{}, while the largest derivative
discrepancy is \JacobianMaximumDerivativeDisagreement{}.  These values satisfy
$\delta_D\le\JacobianDerivativeTolerance{}$ and establish sampled local
admissibility.  They do not imply global injectivity.  The case-wise minima and
representative determinant fields are shown in Fig.~\ref{fig:jacobianaudit}.

\FloatBarrier

Fig.~\ref{fig:jacobianaudit} resolves the distinction between a small
Jacobian and a failed map.  The case-wise minima span several orders of
magnitude because stronger powers deliberately collapse cells toward the tip,
but every sampled value remains positive.  The representative determinant
field likewise decays toward the collapsed edge without changing sign or
forming an interior zero band.  Together with the derivative discrepancy in
Eq.~\eqref{eq:jacobianaudit}, this establishes that the learned density changes
cell concentration without introducing a sampled local inversion.  Only the
planar determinant needed by the reported Galerkin problems is displayed;
auxiliary surface and interface visualizations are omitted because no
corresponding shell or interface weak-form result is claimed.

\paragraph{Cost at a prescribed accuracy.}

\begin{figure}[htp]
\centering
\begin{subfigure}[t]{0.45\linewidth}\centering
\includegraphics[width=\linewidth]{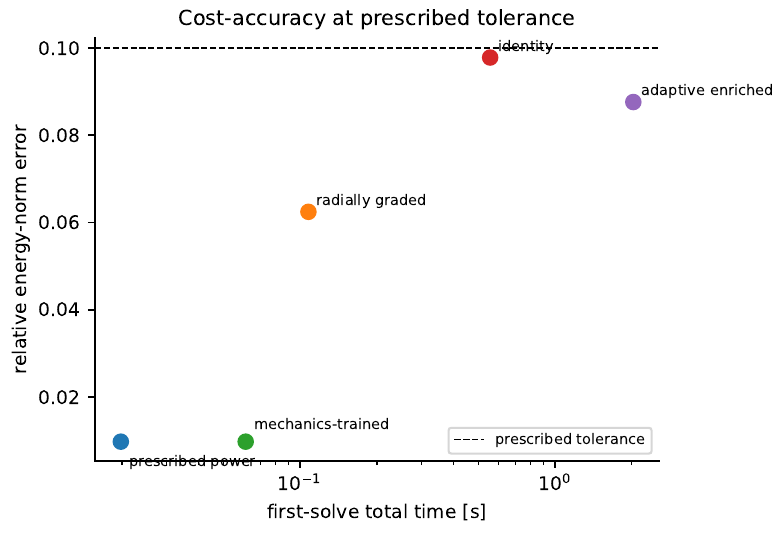}
\caption{First-solve cost versus energy-norm error.}
\end{subfigure}\hfill
\begin{subfigure}[t]{0.45\linewidth}\centering
\includegraphics[width=\linewidth]{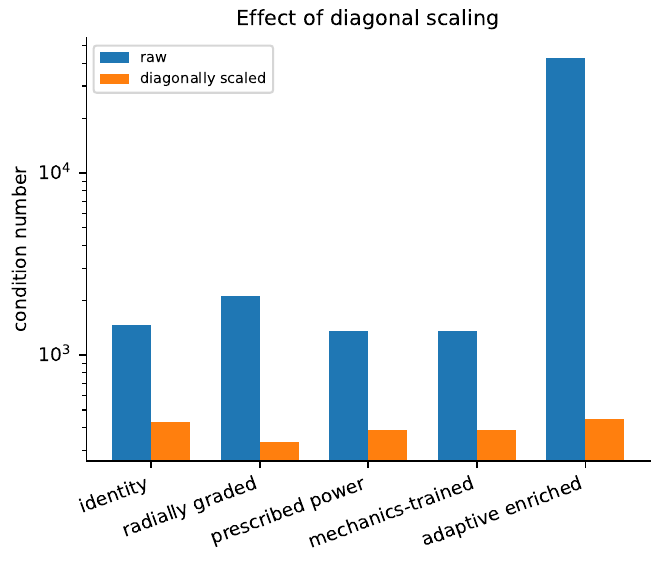}
\caption{CG iterations after diagonal scaling.}
\end{subfigure}

\medskip
\begin{subfigure}[t]{0.45\linewidth}\centering
\includegraphics[width=\linewidth]{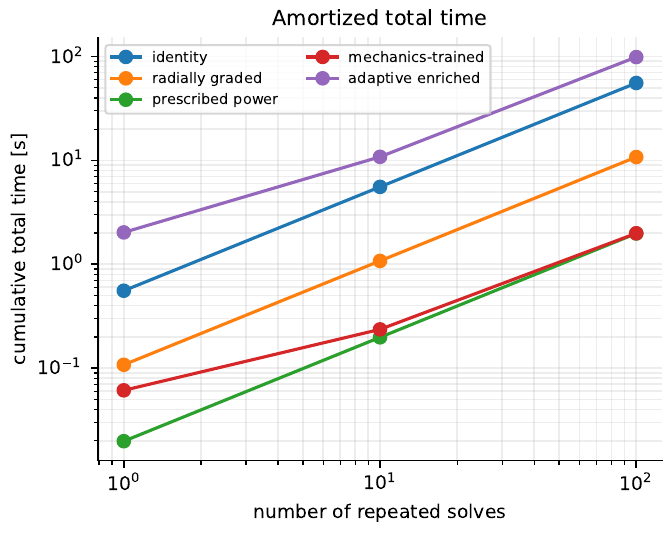}
\caption{Repeated-solve accumulation.}
\end{subfigure}\hfill
\begin{subfigure}[t]{0.45\linewidth}\centering
\includegraphics[width=\linewidth]{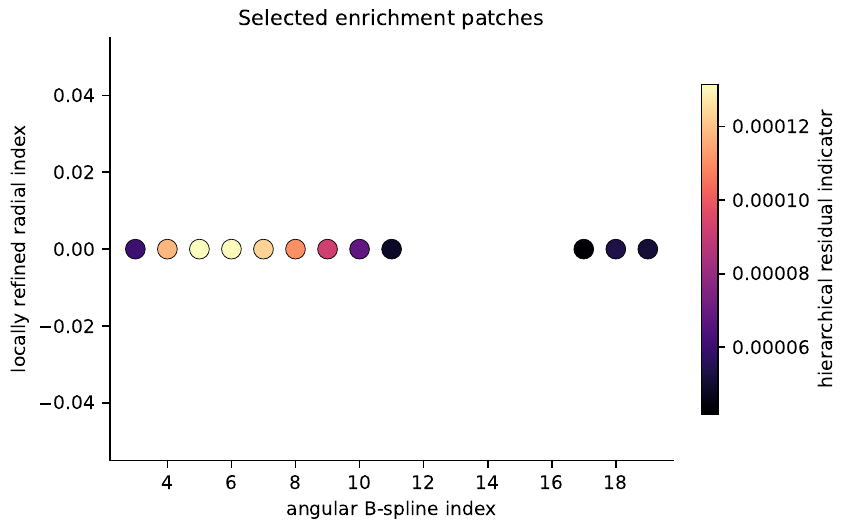}
\caption{Adaptive-patch distribution.}
\end{subfigure}
\caption{Cost--accuracy, diagonal scaling, amortization, and adaptive
enrichment.}
\label{fig:adaptivecost}
\end{figure}

Because coordinate regularization changes the number of degrees of freedom
required at a prescribed accuracy, cost is compared at the prescribed tolerance
$e_E\le\CostTargetEnergy{}$.  Setup and online costs are separated to
distinguish map construction from repeated use.  Wall times are summarized by
their median and interquartile range, and peak working memory is reported for
each computation.  Diagonal scaling reduces the adaptive system's condition
number from \AdaptiveRawCondition{} to \AdaptiveScaledCondition{}, but its
first-solve time remains \AdaptiveFirstSolveSeconds{}~\si{\second} for the
tested enrichment.  Table~\ref{tab:totalcost} reports time in
\si{\second} and memory in \si{\mebibyte}; a dash indicates that no finite
break-even occurs against the identity map under the stated linear
amortization model.  Fig.~\ref{fig:adaptivecost} shows how scaling and the
residual-selected patches contribute to this comparison.

\begin{table}[t]
\centering
\ResultTableFont
\caption{Computational cost, memory, conditioning, CG iterations, and break-even
solve counts at the prescribed energy-norm tolerance.}
\label{tab:totalcost}
\begin{adjustbox}{max width=\linewidth}
\begin{tabular}{lrrrrrrrrr}
\toprule
Method & DOF & $e_E$ & online [s] & setup [s] & first solve [s] & peak [MiB] & raw/scaled $\kappa_2$ & CG iter. & break-even solves \\
\midrule
Identity map & $440$ & $0.0978$ & $0.5254$ & $0$ & $0.5254$ & $311.6$ & $1.470\times 10^{3}$/$4.273\times 10^{2}$ & $159$ & -- \\
Radially graded knots & $256$ & $0.06244$ & $0.1024$ & $0$ & $0.1024$ & $78.98$ & $2.117\times 10^{3}$/$3.349\times 10^{2}$ & $131$ & $1$ \\
Prescribed power map & $120$ & $0.009772$ & $0.01953$ & $0$ & $0.01953$ & $8.998$ & $1.358\times 10^{3}$/$3.88\times 10^{2}$ & $102$ & $1$ \\
Mechanics-trained power map & $120$ & $0.009772$ & $0.01968$ & $0.04181$ & $0.06149$ & $8.997$ & $1.358\times 10^{3}$/$3.88\times 10^{2}$ & $102$ & $1$ \\
Adaptive enriched B-spline & $464$ & $0.0876$ & $0.9194$ & $0.9951$ & $1.914$ & $836.1$ & $4.302\times 10^{4}$/$4.467\times 10^{2}$ & $164$ & -- \\
\bottomrule
\end{tabular}
\end{adjustbox}
\end{table}

Once trained, the mechanics-trained power map has essentially the same online
cost as the prescribed power map; its additional cost is confined to setup.  Repeated use
therefore reduces the relative contribution of training, whereas the adaptive
comparator retains both a larger setup and a larger online cost at the tested
target.  Table~\ref{tab:amortization} accumulates these contributions and gives
\AdaptiveBreakEven{} for the adaptive break-even result.

\begin{table}[htp]
\centering
\ResultTableFont
\caption{Amortized total solution time.}
\label{tab:amortization}
\begin{adjustbox}{max width=\linewidth}
\begin{tabular}{lrrr}
\toprule
Method & 1 solve (\si{\second}) & 10 solves (\si{\second}) & 100 solves (\si{\second}) \\
\midrule
Identity map & $0.5254$ & $5.254$ & $52.54$ \\
Radially graded knots & $0.1024$ & $1.024$ & $10.24$ \\
Prescribed power map & $0.01953$ & $0.1953$ & $1.953$ \\
Mechanics-trained power map & $0.06149$ & $0.2386$ & $2.01$ \\
Adaptive enriched B-spline & $1.914$ & $10.19$ & $92.94$ \\
\bottomrule
\end{tabular}
\end{adjustbox}
\end{table}

\FloatBarrier

Fig.~\ref{fig:adaptivecost} converts the accuracy comparison into a
computational comparison at the prescribed error tolerance.  The first panel includes both
setup and online work needed to reach the target error, while the second shows
that diagonal scaling reduces iterative difficulty but does not remove the
adaptive method's larger CG iteration count at this target.  The repeated-solve curves expose
the different intercepts and slopes: training is a one-time cost for the
mechanics-trained map, whereas its online cost remains close to the prescribed-power
solve.  The adaptive-patch distribution explains the larger memory and system
size by showing where enrichment is retained.  This figure is important
because equal-DOF accuracy alone cannot establish efficiency; the conclusion
must persist when every method is compared at the same error requirement.

\paragraph{Assembly on the collapsed-edge trial space.}

Direct assembly on $V_h^{\rm adm}$ excludes collapsed-edge basis functions
before matrix construction.  It agrees with constrained full assembly to
$1.332\times10^{-15}$ in the vector problem and gives no stored entries associated with the excluded basis functions in the scalar and adaptive systems.  The collapsed-edge treatment is
therefore algebraically consistent across the three discretizations.
Table~\ref{tab:tipassembly} reports the excluded basis counts and the
corresponding maximum coefficient or entry difference.

\begin{table}[htp]
\centering
\ResultTableFont
\caption{Collapsed-edge trial-space assembly checks.}
\label{tab:tipassembly}
\begin{adjustbox}{max width=\linewidth}
\begin{tabular}{lrrlr}
\toprule
System & excluded collapsed-edge bases & post-assembly row deletions & verification quantity & value \\
\midrule
Scalar & $21$ & 0 & max. stored excluded-row entry & $0$ \\
Vector & $26$ & 0 & max. coefficient difference & $1.332\times 10^{-15}$ \\
Adaptive & $44$ & 0 & max. stored excluded-row entry & $0$ \\
\bottomrule
\end{tabular}
\end{adjustbox}
\end{table}


\section{Conclusions}
\label{sec:conclusions}

Restricting the network to a positive radial density yields a singular
coordinate that is radially monotone by construction and trainable through
Galerkin functionals evaluated at discrete equilibrium.  The crack calculations show that this
coordinate transfers consistently to elasticity and preserves the stress intensity factors obtained from the annulus-based displacement fit and
the interaction integral.  In the
nonlinear Robin family, the density correction remains nonzero under two-sided
extrapolation and reduces the field, energy, and near-tip amplitude errors for
every accepted initialization.

The benefit depends on whether the required coordinate is already known.  The
straight crack is resolved by selecting the classical square map, and the
affine wedge is represented more accurately by the affine grading-exponent predictor.
Neural flexibility becomes useful when the coordinate contains a non-power-law
correction that cannot be prescribed from the leading exponent alone.  The
quadrature, Jacobian, collapsed-edge trial-space assembly, and prescribed-accuracy tests show
that this distinction persists beyond the training functional.

The present conclusions are limited to manufactured quasi-static problems and
to a straight-crack, homogeneous-material interaction integral.  Extending the method to an
evolving, curved, interfacial, or materially nonlinear crack requires training
against the corresponding weak form and fracture functional; it is not implied
by the coordinate tests reported here.  Such extensions address a broader
problem than neural solvers based on enriched local fields, peridynamic energy,
phase-field constitutive models, or finite-element-integrated networks
\cite{Gu2023EnrichedPINN,Zhu2025ExtremeFracture,Ning2023PeridynamicPINN,
Dammass2025Hybrid,Pantidis2026IFENN}.

\bibliographystyle{unsrt}
\bibliography{references}

\end{document}